\documentclass[11pt,reqno]{amsart}
\usepackage[T1]{fontenc}
\usepackage{lmodern,amssymb,microtype}
\usepackage[letterpaper,left=1.08in,right=1.08in,top=.96in,bottom=.98in,
 headheight=14pt,headsep=15pt]{geometry}
\newcommand{\amend}[1]{#1}
\usepackage[hidelinks]{hyperref}
\hypersetup{pdftitle={Three problems of Farmer and Johnson on Lipschitz summability},
 pdfauthor={Daniel M. Pellegrino and Eduardo V. Teixeira},
 pdfkeywords={Lipschitz summing maps, Farmer-Johnson problems, Grothendieck theorem, Lipschitz extension, Hamming cube}}
\newtheorem{theorem}{Theorem}[section]
\newtheorem{lemma}[theorem]{Lemma}
\newtheorem{proposition}[theorem]{Proposition}
\newtheorem{corollary}[theorem]{Corollary}
\DeclareMathOperator{\Lip}{Lip}
\DeclareMathOperator{\diam}{diam}
\DeclareMathOperator{\Var}{Var}
\DeclareMathOperator{\supp}{supp}
\newcommand{\R}{\mathbb R}
\newcommand{\Z}{\mathbb Z}
\newcommand{\E}{\mathbb E}
\newcommand{\pL}{\pi_p^{\mathrm L}}
\newcommand{\Id}{\mathrm{Id}}
\newcommand{\norm}[1]{\lVert#1\rVert}

\newcommand{\cD}{\mathcal D}
\newcommand{\dH}{d_{\mathrm H}}
\allowdisplaybreaks[1]
\title[Three problems on Lipschitz summability]{Three problems of Farmer and Johnson\\ on Lipschitz summability}
\author{Daniel M. Pellegrino}
\address{Departamento de Matem\'atica, Universidade Federal da Para\'iba, Jo\~ao Pessoa, PB, Brazil}
\author{Eduardo V. Teixeira}
\address{Department of Mathematics, Oklahoma State University, Stillwater, Oklahoma, USA}
\date{}
\numberwithin{equation}{section}
\begin{document}
\begin{abstract}
Farmer and Johnson asked which properties of absolutely summing
operators persist for Lipschitz maps. We give negative answers to their
Problems~2, 4, and~5. The Grothendieck assertions for maps from $L_1$
or $C(K)$ into Hilbert space fail even for smooth Lipschitz maps with
relatively compact range. Lipschitz $2$-summability does not imply
Lipschitz $2$-integrality, even for a compact countable domain and a
reflexive target of cotype two. A six-point metric space disproves the
classical conjugate-exponent norm test for $1<p<\infty$.
\amend{Finite subsets of $\ell_1$ recover the Lipschitz $p$-summing norm
on every metric domain, and ordinary Hamming cubes recover the absolutely
$p$-summing norm of every bounded linear operator.}
\end{abstract}
\maketitle

\section{Introduction}

The theory of absolutely summing operators grew out of
Grothendieck's work on tensor products \cite{Groth} and Pietsch's
domination and factorization theory \cite{Pietsch}; see also
\cite{DJT}. Grothendieck's theorem implies that every bounded linear
operator from $L_1(\mu)$ into a Hilbert space is absolutely
$1$-summing, and that every bounded linear operator from $C(K)$, where $K$ is compact Hausdorff, into
a Hilbert space is absolutely $2$-summing. The norm of an absolutely
$p$-summing linear operator also admits a test by precomposition
with operators defined on $\ell_{p'}$, where $1/p+1/p'=1$
when $p>1$, and $p'=\infty$ when $p=1$.

Farmer and Johnson \cite{FJ} extended $p$-summability to Lipschitz
maps between metric spaces. In their definition, pairs of points
take the place of vectors, and scalar Lipschitz functions replace
linear functionals. They established a nonlinear Pietsch domination
and factorization theorem and proved that, for linear operators, the
new norm coincides with the classical one. They also obtained a
Grothendieck-type theorem for metric trees. The nonlinear factorization does not automatically provide an
extension of the factored map to the whole intermediate space.

Lipschitz-free spaces furnish a linear setting for the study of
Lipschitz maps; foundational results are due to Godefroy and Kalton
\cite{GK}, and Godard \cite{Godard} identified their connection
with tree metrics. Further aspects of Lipschitz summability include
the extrapolation results of Chen and Zheng \cite{CZ}, the duality
theory of Ch\'avez-Dom\'inguez \cite{CD}, and the work on strongly
summing maps by Saadi \cite{Saadi}. Saleh \cite{Saleh} computed
summing norms in finite metric examples. Masty\l{}o and
S\'anchez P\'erez \cite{MS} obtained domination theorems for
certain Lipschitz maps on $C(K)$ under additional hypotheses.
For Lipschitz extension into Banach spaces, Ball's theorem
\cite{Ball} and the work of Mendel and Naor \cite{MN} provide
estimates under geometric hypotheses on the range.

Problems~2, 4, and~5 in \cite{FJ} ask, respectively, whether
Lipschitz $2$-summability implies Lipschitz $2$-integrality;
whether the two Hilbert-valued Grothendieck coincidences persist
for arbitrary Lipschitz maps on $L_1$ and $C(K)$; and whether the
Lipschitz $p$-summing norm admits the classical
$\ell_{p'}$-precomposition test. The first implication fails even
for a compact countable subset of $c_0$ with values in a reflexive
space of Rademacher cotype two. Both Hilbert-valued assertions
fail, already for smooth contractions with relatively compact
range. For $1<p<\infty$, a six-point metric space rules out the
proposed $\ell_{p'}$ test.

There are nevertheless exact norm-testing formulas. For maps on
Banach spaces, contractions defined on finite-dimensional sup-norm
spaces recover the Lipschitz $p$-summing norm. For arbitrary metric
domains, finite subsets of $\ell_1$ suffice. The absolutely
$p$-summing norm of a linear operator is recovered by composition
with contractions on $\ell_1$ or on finite Hamming cubes. These
formulas also yield quantitative estimates for Hilbert-valued maps.

Let $(X,d_X)$ and $(Y,d_Y)$ be metric spaces, and let $T:X\to Y$ be
Lipschitz. For a scalar-valued map $g$ on $X$, let $\Lip(g)$ denote
its Lipschitz constant and put
$\Lip_1(X)=\{g:X\to\R:\Lip(g)\le1\}$.
For $1\le p<\infty$, the Lipschitz $p$-summing norm $\pL(T)$
is the least $C\in[0,\infty]$ for which
\begin{equation}\label{eq:lippsum}
 \left(\sum_{i=1}^N w_i d_Y(Tx_i,Ty_i)^p\right)^{1/p}
 \le C\sup_{g\in\Lip_1(X)}
       \left(\sum_{i=1}^N w_i |g(x_i)-g(y_i)|^p\right)^{1/p}
\end{equation}
holds for every $N\ge1$, $x_i,y_i\in X$, and $w_i\ge0$.
Real-valued contractions can be normalized at any fixed point by subtracting a
constant, and a contraction defined on a subset extends to $X$ with
the same constant \cite{McShane}. Consequently, restriction to an
isometric subset preserves the scalar supremum for tests supported
there. We write $\cD_{p,X}((x_i,y_i,w_i)_i)$ for the supremum on the
right-hand side of \eqref{eq:lippsum}, and $\E$ for expectation.
All Banach spaces and scalar Lipschitz functions are real. For pointed
maps, a distinguished point is mapped to the distinguished point of
the target. The conjugate exponent is given by $1/p+1/p'=1$ for
$p>1$, with $p'=\infty$ when $p=1$. The Hamming cube
$Q_N=\{0,1\}^N$ has distance $\dH(u,v)=|\{j:u_j\ne v_j\}|$.
For a bounded linear operator $U:X\to Y$, write $B_{X^*}$ for
the closed unit ball of $X^*$. The classical absolutely $p$-summing
norm $\pi_p(U)$ is the least $C$ such that
\begin{equation}\label{eq:linearpsum}
 \left(\sum_{i=1}^N\|Uu_i\|^p\right)^{1/p}
 \le C\sup_{x^*\in B_{X^*}}\left(\sum_{i=1}^N|x^*(u_i)|^p\right)^{1/p}
\end{equation}
for every finite family $u_1,\ldots,u_N\in X$.
For bounded linear maps, $\pL(U)=\pi_p(U)$ by
\cite[Theorem~2]{FJ}. The definition of Lipschitz $p$-integrality,
used in Theorem~\ref{thm:integral}, is given in
Section~\ref{sec:integral}.

\begin{theorem}\label{thm:testing}
Let $X$ be a Banach space, $Y$ a metric space, and $T:X\to Y$ a Lipschitz
map. For every $1\le p<\infty$,
\begin{equation}\label{eq:global-test}
 \pL(T)=\sup_{n\ge1}\ \sup_{\substack{F:\ell_\infty^n\to X\\ \Lip(F)\le1}}
                   \pL(T\circ F).
\end{equation}
\end{theorem}

In \eqref{eq:global-test}, each $F$ is defined on the whole space
$\ell_\infty^n$, not merely on a finite subset. The isometric
embedding of a finite metric space into $\ell_\infty^n$ does not,
by itself, give a Lipschitz extension of arbitrary Banach-valued
maps. The extensions needed for \eqref{eq:global-test} are obtained
from finite graphs, with Lipschitz constants tending to one.

\begin{theorem}\label{thm:ck}
Let $K$ be an infinite compact Hausdorff space. In
Theorem~\ref{thm:testing} the same supremum is obtained by allowing
contractions $F:C(K)\to X$. For every infinite-dimensional Banach
space $Y$ there is a Fr\'echet $C^\infty$ contraction $U:C(K)\to Y$
with relatively compact range such that
\[
 \pL(U)=\infty\qquad(1\le p<\infty).
\]
\end{theorem}

Already for a real Hilbert target of dimension $d$, the optimal
$2$-summing bound for Lipschitz contractions on $C(K)$ is $\sqrt d$.
The linear bound is independent of $d$. In a different direction,
Masty\l o and S\'anchez P\'erez \cite{MS} prove domination results
for Lipschitz maps on $C(K)$ under concavity and summability hypotheses.
Theorem~\ref{thm:ck} concerns all Lipschitz maps.

\begin{theorem}\label{lone:thm:failure}
For every infinite-dimensional Banach space $Y$ there exists a
Fr\'echet $C^\infty$ contraction $T:\ell_1\to Y$ with relatively
compact range such that
\begin{equation}\label{eq:lone-all-p}
 \pL(T)=\infty\qquad(1\le p<\infty).
\end{equation}
In particular, the $L_1$-to-Hilbert question in
\cite[Problem~4]{FJ} has a negative answer already for
$\ell_1\to\ell_2$.
\end{theorem}

\begin{theorem}\label{lone:thm:recovery}
Let $U:X\to Y$ be a bounded linear operator between Banach spaces
and let $1\le p<\infty$. Then
\begin{equation}\label{lone:eq:recovery}
 \pi_p(U)
 =\max_{\substack{F:\ell_1\to X\\ \Lip(F)\le1}}\pL(U\circ F)
 =\sup_{N\ge1}\ \sup_{\substack{F:Q_N\to X\\ \Lip(F)\le1}}
       \pL(U\circ F).
\end{equation}
The maximum exists also when $\pi_p(U)=\infty$.
\end{theorem}

\amend{Taking $U=\Id_{\ell_2^d}$ in Theorem~\ref{lone:thm:recovery}
and using Corollary~\ref{lone:cor:rank} gives
$$
 \max_{\substack{F:\ell_1\to\ell_2^d\\ \Lip(F)\le1}}\pi_1^{\mathrm L}(F)
 =\frac{\sqrt{\pi}\,\Gamma((d+1)/2)}{\Gamma(d/2)}
 \sim\sqrt{\frac{\pi d}{2}}\qquad(d\to\infty).
$$
For linear contractions $F$ with the same source and target,
Grothendieck's theorem gives $\pi_1(F)\le K_G$ in every dimension,
where $K_G$ is the real Grothendieck constant.}

Chen and Zheng \cite{CZ} obtained summability results under
additional assumptions on the Lipschitz adjoint. The arXiv version
of Saadi \cite[Corollary~3.6]{Saadi} includes an affirmative assertion
for maps from $\ell_1$ into Hilbert space, based on an isometric
identification of the Lipschitz-free space $\mathcal F(\ell_1)$ with
an $L_1$ space \amend{asserted in} \cite[Corollary~8]{DTZ}.
\amend{For the Hamming square $\{0,e_1,e_2,e_1+e_2\}\subset\ell_1$,
the three ways of partitioning the vertices into two pairs give sums
of distances $4,2,2$. Godard's criterion \cite[Theorem~4.2]{Godard}
requires the largest two sums to agree, which rules out the
isometric identification.}

\begin{theorem}\label{thm:integral}
There exist a compact countable set $M\subset c_0$, a separable
reflexive Banach space $Y$ of Rademacher cotype two, and a Lipschitz
$1$-summing map $T:M\to Y$ that has no Lipschitz extension from $c_0$
into $Y^{**}$. In particular, $T$ is Lipschitz $2$-summing and fails
to be Lipschitz $q$-integral for every $1\le q<\infty$.
\end{theorem}

In the Lipschitz Pietsch factorization the final map is initially
defined on a subset of $L_2$. A Hilbert target permits its extension
to all of $L_2$ by Kirszbraun's theorem. Farmer and Johnson also give
a positive result for targets with Ball's Markov cotype two property
\cite{FJ,Ball}. Theorem~\ref{thm:integral} shows that Rademacher
cotype two does not suffice. Its finite examples are obtained from
nets of a Euclidean sphere. Two explicit families of real-valued
contractions control their summing norms. Extending the vector-valued
map would force a trace estimate incompatible with Grothendieck's
inequality and the divergence theorem.

\begin{theorem}\label{thm:six}
Let $M$ be the six vertices of $K_{3,3}$ with the graph distance. For
all $1\le p<\infty$ one has $\pL(\Id_M)=3^{1/p}$. For each such
$p$ and each fixed uniformly convex Banach space $V$,
\[
 \sup_{\substack{A\subset V\ \mathrm{finite}\\ S:A\to M,\ \Lip(S)\le1}}
       \pL(S)<3^{1/p}.
\]
Consequently, the test with $V=\ell_{p'}$ fails for
$1<p<\infty$.
\end{theorem}

For a uniformly convex space, three points with mutual distances
at least two cannot all lie in a ball of radius $1+\varepsilon$
when $\varepsilon>0$ is sufficiently small. This yields a uniform
strict inequality for every pullback of $K_{3,3}$. At $p=1$, the
original test with $\ell_\infty$ holds. Finite subsets of $\ell_1$
give an exact test for all $1\le p<\infty$, and in the case of Banach
domains the supremum is attained by a map defined on the whole of
$\ell_1$.

\section{Testing on finite-dimensional sup-norm spaces}\label{ck:section}

Fix a finite metric space $S$. The vertices of the graph will be
partitioned into equally sized fibers indexed by $S$, with matchings
realizing the distances between labels. Expansion within each fiber
controls scalar oscillation. A graph has edge expansion at least $h$
if every nonempty vertex set of cardinality at most half the total
has at least $h$ times its cardinality in boundary edges.

\begin{lemma}\label{ck:fibers}
Let $S$ be a finite metric space with at least two points. For every
$\delta>0$ and $D>0$ there is a finite connected metric graph $G$ of
girth greater than $D$, whose vertex set $V$ is partitioned into sets
$V_x$, $x\in S$, of equal size. For every two distinct
$x,y\in S$, there is a perfect matching between $V_x$ and $V_y$
consisting of edges of length $d_S(x,y)$.
Within each $V_x$ there is a cubic graph with edges of length
$\delta$ and edge expansion bounded below by an absolute $h>0$.
\end{lemma}

\begin{proof}
Write $s=|S|$ and take the Cartesian product of graphs $K_s\square K_4$.
In a random $n$-lift each base vertex is replaced by $n$ vertices and
independent uniform perfect matchings replace its edges. Above each
$x\in S$, the four base vertices thus give a random $n$-lift of $K_4$.
Since $K_4$ has six edges and four vertices, the expansion theorem
of Amit and Linial \cite[Theorem~1.1]{AL} applies to each of its
random lifts. Thus, for a fixed $h>0$ depending only on $K_4$,
the probability that all $s$ vertical subgraphs have edge
expansion at least $h$ tends to one. The base graph
$K_s\square K_4$ is fixed, connected, and loop-free. By
\cite[Lemma~5.1]{GJR}, the numbers $Z_k$ of cycles of lengths
$2\le k\le m$ in its random lift converge jointly to independent
Poisson random variables of finite means $\lambda_k$.
Consequently,
\[
 \lim_{n\to\infty}\mathbb P(Z_2=\cdots=Z_m=0)
 =\exp\Big(-\sum_{k=2}^{m}\lambda_k\Big)>0.
\]
For fixed $m$, this positive limiting probability and the
asymptotically certain expansion event have a nonempty intersection
for all sufficiently large $n$. The vertical graphs are
connected by expansion. Between any two labels the four horizontal
matchings together give a perfect matching, so the whole graph is
connected. For each horizontal edge connecting distinct labels $x,y$ assign
length $d_S(x,y)$, and assign length $\delta$ to each vertical edge.
Choose $m$ so that
\begin{equation}\label{eq:girth-choice}
 m\min\bigl\{\delta,\min_{x\ne y}d_S(x,y)\bigr\}>D.
\end{equation}
Realize the edges by intervals with the induced path metric. Every
simple cycle has at least $m+1$ edges, hence length greater than $D$
by \eqref{eq:girth-choice}.
\end{proof}

Let $\mu_x$ denote uniform probability on $V_x$, and set
$\bar g(x)=\int_{V_x}g\,d\mu_x$ for a scalar contraction
$g:V\to\R$. Fix $x\ne y$, and use the matching to identify
$V_x$ with $V_y$. Integrating the pointwise inequality along the
matching yields
\begin{equation}\label{eq:fiber-means}
 |\bar g(x)-\bar g(y)|\le d_S(x,y).
\end{equation}
We also have, for each $1\le p<\infty$,
\begin{equation}\label{eq:fiber-oscillation}
 \|g-\bar g(x)\|_{L_p(\mu_x)}\le C_p\delta,
\end{equation}
where $C_p$ depends only on $p$ and the absolute expansion constant
$h$, not on $S$, the size of the fibers, or $\delta$.

Put $q=1+h/3>1$.
Every set $A$ containing at most half the vertices of a vertical
cubic graph has a closed one-step neighborhood of cardinality at
least $q|A|$: there are at least $h|A|$ boundary edges, and a
vertex outside $A$ meets at most three such edges. Choose a median
$a$ of $g|_{V_x}$. Since $g$ changes by at most $\delta$ along each
vertical edge, iteration of the neighborhood inequality gives
\begin{equation}\label{eq:fiber-tail}
 \mu_x\bigl(\{|g-a|>k\delta\}\bigr)\le q^{-k}
 \qquad(k=0,1,\ldots).
\end{equation}
Indeed, let
$A_k^+=\{v:g(v)>a+k\delta\}$ and
$A_k^-=\{v:g(v)<a-k\delta\}$. Both $A_0^+$ and $A_0^-$
have at most half the vertices. Every neighbor of a vertex in
$A_{k+1}^+$ belongs to $A_k^+$, and the same assertion holds for
the negative tails. Expansion therefore gives
$|A_k^\pm|\ge q|A_{k+1}^\pm|$. Induction yields
$\mu_x(A_k^\pm)\le (2q^k)^{-1}$, and summing the two estimates
proves \eqref{eq:fiber-tail}.
By layer-cake integration,
\begin{equation}\label{eq:fiber-moment}
 \int_{V_x}|g-a|^p\,d\mu_x
 \le\delta^p\sum_{k=0}^\infty
       \bigl((k+1)^p-k^p\bigr)q^{-k}<\infty.
\end{equation}
Finally, $|\bar g(x)-a|\le\|g-a\|_{L_p(\mu_x)}$, and the
triangle inequality in $L_p(\mu_x)$ yields
\eqref{eq:fiber-oscillation}. In particular, one may take
\[
 C_p=2\left[\sum_{k=0}^{\infty}
 \big((k+1)^p-k^p\big)q^{-k}\right]^{1/p}.
\]

\begin{lemma}\label{ck:tube}
Given $L>1$ and $\rho>0$, every compact finite metric graph $G$ of
sufficiently large girth has the following property. In any metric
space $E$ containing $G$ isometrically, the closed
$\rho$-neighborhood of $G$ admits an $L$-Lipschitz retraction onto $G$.
\end{lemma}

\begin{proof}
Set
\begin{equation}\label{eq:tube-parameters}
 M=\frac{2(L+1)}{L-1},\qquad a=L(M+1)\rho.
\end{equation}
It suffices that the metric girth of $G$ exceed $4a$. Every closed ball $B_G(b,a)$ is then a compact metric tree with
the metric inherited from $G$. Two distinct shortest paths from $b$
to the same point of the ball would contain a cycle of length at most
$2a$, so radial geodesics are unique. Two such geodesics cannot split
and subsequently meet, and their union is a rooted tree. A path
between two points of the ball through their last common ancestor
has length at most $2a$. Any strictly shorter path in $G$, or any
different path of equal length, would produce a cycle of length at
most $4a$. Consequently the tree paths are shortest paths in $G$.

Start with the identity on $G$ and consider a partial $L$-Lipschitz
extension $r$ from $D$ into $G$, where $D$ lies in the tube and
contains $G$. Given a new point $u$ in the tube, choose a nearest
point $b\in G$. For $v\in D$ with $d_E(u,v)\le M\rho$, the points
$r(v)$ lie in the tree $B_G(b,a)$, because
\begin{equation}\label{eq:tube-near}
 d_G(r(v),b)\le Ld_E(v,b)
 \le L\bigl(d_E(v,u)+d_E(u,b)\bigr)\le a.
\end{equation} In this tree form the closed
balls with centers $r(v)$ and radii $L d_E(u,v)$. Any two intersect
by the triangle inequality and the Lipschitz bound on $D$.
The balls intersect pairwise within the tree. Indeed, for two centers
$r(v),r(w)$ their distance is at most
$L(d_E(u,v)+d_E(u,w))$. Since the tree is isometrically embedded,
the geodesic joining these centers is contained in the tree; the
distance inequality provides a point of that geodesic in both
balls. This also covers radii exceeding the diameter of the tree. Closed balls in a metric tree are closed subtrees, and
pairwise intersecting closed subtrees of a compact tree have a
common point. For a finite family, root the tree and choose a
subtree whose point nearest to the root is farthest away. The
geodesic from the root to a point in its intersection with any
other subtree contains that chosen point. Compactness gives a
common point for the whole family.

Choose a common point $z$. The constraint at $v=b$ gives
$d_G(z,b)\le L\rho$. For any remaining $v\in D$, choose a nearest point $c\in G$
to $v$. Since $d_G(r(v),c)\le L\rho$ and
\[
 d_G(b,c)=d_E(b,c)
 \le d_E(b,u)+d_E(u,v)+d_E(v,c)
 \le d_E(u,v)+2\rho,
\]
we have
\begin{align}
 d_G(z,r(v))
 &\le d_G(z,b)+d_G(b,c)+d_G(c,r(v))\notag\\
 &\le d_E(u,v)+2(L+1)\rho\notag\\
 &\le Ld_E(u,v).\label{eq:tube-far}
\end{align}
The last inequality follows from $d_E(u,v)>M\rho$ and
\eqref{eq:tube-parameters}. Setting
$r(u)=z$ preserves the Lipschitz bound. The union of a chain of these
partial extensions is again such an extension, so the maximality
principle completes the retraction.
\end{proof}

The condition $L>1$ cannot in general be replaced by $L=1$. Miesch
exhibited an isometric embedding of an intrinsic circle for which
every neighborhood fails to admit a $1$-Lipschitz retraction
\cite[Example~4.2]{Miesch}.

\begin{proof}[Proof of Theorem~\ref{thm:testing}]
Precomposition with a contraction can only decrease $\pL$. To prove
the reverse inequality, take a finite weighted test
$(x_i,y_i,w_i)_i$ for $T:X\to Y$, and let $S\subset X$ be its set
of endpoints together with $0_X$. We may assume $|S|\ge2$. Write
\begin{equation}\label{eq:test-numerators}
 A=\left(\sum_i w_i d_Y(Tx_i,Ty_i)^p\right)^{1/p},\qquad
 B=\sup_{\Lip(f)\le1}
     \left(\sum_i w_i|f(x_i)-f(y_i)|^p\right)^{1/p},
\end{equation}
where $f:S\to\R$. This is also the scalar supremum on $X$, by
real-valued Lipschitz extension. Fix $\eta,\delta>0$ and choose
$\rho$ so large that $\max_{x\in S}\|x\|/\rho\le\eta$. Apply
Lemma~\ref{ck:fibers} with the girth required by Lemma~\ref{ck:tube}
for this $\rho$ and $L=1+\eta$.

Map $V_x$ to $x$ for every $x\in S$ and interpolate affinely on
each edge of $G$. The resulting $b:G\to X$ is a contraction and
$\|b\|_\infty\le\max_{x\in S}\|x\|$. Embed $G$ isometrically in
$E=\ell_\infty(G)$ by $j(v)=d_G(v,\cdot)$. On its closed tube
use the retraction $r$ of Lemma~\ref{ck:tube}, identifying $G$
with $j(G)$. The map
\begin{equation}\label{eq:cutoff-extension}
 \widetilde b(z)=
 \begin{cases}
 \bigl(1-d_E(z,G)/\rho\bigr)b(r(z)),&d_E(z,G)\le\rho,\\
 0,&d_E(z,G)>\rho
 \end{cases}
\end{equation}
is defined on all of $E$. On the tube its Lipschitz
constant is at most
\begin{equation}\label{eq:cutoff-bound}
 \Lip(\widetilde b)\le L+\frac{\|b\|_\infty}{\rho}
 \le1+2\eta,\qquad L=1+\eta.
\end{equation}
The same estimate holds when one point, say $z'$, is outside
the tube. In that case $\widetilde b(z')=0$ and, writing
$t=d_E(z,G)\le\rho$, we obtain
\[
 \|\widetilde b(z)-\widetilde b(z')\|
 \le \|b\|_\infty\frac{\rho-t}{\rho}
 \le\frac{\|b\|_\infty}{\rho}d_E(z,z'),
\]
because $d_E(z',G)>\rho$ and the distance to $G$ is
$1$-Lipschitz.

To pass to a finite-dimensional source, choose $v_0\in V_{0_X}$
and set $\iota(v)=(d_G(v,a)-d_G(v_0,a))_{a\in V}$. This embeds
$V$ isometrically in $\ell_\infty^V$ and sends $v_0$ to zero.
The map $\iota(v)\mapsto j(v)$ extends to a contraction
$P:\ell_\infty^V\to E$. Explicitly, for $z\in\ell_\infty^V$
and $s\in G$ set
\begin{equation}\label{eq:coordinate-extension}
 P(z)(s)=\min_{v\in V}\left(d_G(v,s)+\|z-\iota(v)\|_\infty\right).
\end{equation}
For each fixed $s$ this is a scalar $1$-Lipschitz extension
of $\iota(v)\mapsto d_G(v,s)$: at $\iota(v)$ the term indexed by
$v$ gives $d_G(v,s)$, and every other term is no smaller by
the triangle inequality. Hence
$P(\iota(v))=j(v)$, and taking the supremum over $s$ shows
$\|P(z)-P(z')\|_\infty\le\|z-z'\|_\infty$.
In particular, since $\iota(v_0)=0$, the triangle inequality gives
\[
 P(0)(s)=\min_{v\in V}\big(d_G(v,s)+d_G(v,v_0)\big)
       =d_G(v_0,s)=j(v_0)(s).
\]
Each $P(z)$ is bounded since $G$ is compact and $V$ is
finite. Composing $P$ with the preceding extension gives
$F_0:\ell_\infty^V\to X$ with
$F_0(\iota(v))=x$ on $V_x$ and $F_0(0)=0_X$. Set $c=1+2\eta$ and $F(z)=F_0(z/c)$.
The map $F$ is a contraction on the entire sup-norm space, fixes
zero, and has the original labels as its values at $c\,\iota(v)$.

For each $i$ choose $U_i,V_i$ independently with respective laws
$\mu_{x_i},\mu_{y_i}$. Test $T\circ F$ on the pairs
$(c\,\iota(U_i),c\,\iota(V_i))$, weighting each possible pair by
its probability times $w_i$. Its vector numerator is $A$. For a scalar
contraction $h$ on $\ell_\infty^V$ set $g(v)=h(c\,\iota(v))/c$.
This is a contraction on the graph vertices. Its means are a scalar
contraction on $S$, and the fiber estimate controls both deviations. Minkowski's inequality and \eqref{eq:fiber-oscillation} give
\begin{equation}\label{eq:lift-test}
 \left(\sum_i w_i\E|h(c\,\iota(U_i))-h(c\,\iota(V_i))|^p\right)^{1/p}
 \le c\left(B+2C_p\delta\Big(\sum_i w_i\Big)^{1/p}\right).
\end{equation}
If the right-hand side of \eqref{eq:global-test} is finite and
equal to $K$, the definition \eqref{eq:lippsum} and
\eqref{eq:lift-test} imply
\begin{equation}\label{eq:test-conclusion}
 A\le K(1+2\eta)
 \left(B+2C_p\delta\Big(\sum_iw_i\Big)^{1/p}\right).
\end{equation}
The graph and the test map were constructed for arbitrary
$\eta,\delta>0$. Letting both tend to zero yields $A\le KB$
for every finite test. An infinite $K$ is covered by
the first inequality.
\end{proof}

\section{Smooth maps on \texorpdfstring{$C(K)$}{C(K)}}\label{sec:consequences}

Fix an infinite compact Hausdorff space $K$ and a point $t_*\in K$.
Evaluation at distinct points $s_1,\ldots,s_n$ gives a linear
contraction $P:C(K)\to\ell_\infty^n$ with an isometric right inverse.
Indeed, choose continuous functions $0\le\psi_i\le1$ with disjoint
supports and $\psi_i(s_i)=1$. Such functions exist by normality,
using pairwise disjoint neighborhoods of the finite set of points.
Then $z\mapsto\sum_i z_i\psi_i$ is the required section. It follows
that $\pL(U\circ P)=\pL(U)$ for every Lipschitz map $U$ on
$\ell_\infty^n$, since precomposition by either $P$ or its section
has norm at most one. This proves the first assertion of
Theorem~\ref{thm:ck}.

Evaluation at $t_*$ separates the bounded smooth maps used in the
next lemma.

\begin{lemma}\label{lem:glue}
Let $H_j:\ell_\infty^{n_j}\to X$ be bounded $C^\infty$ contractions
into a Banach space, and let $b_j>1$. There are a $C^\infty$
contraction $F:C(K)\to X$ vanishing at zero and affine isometric
embeddings $J_j:\ell_\infty^{n_j}\to C(K)$ such that
\[
 F(J_j(b_j z))=H_j(z)\qquad(z\in\ell_\infty^{n_j}).
\]
If each $H_j$ has relatively compact range and
$\|H_j\|_\infty\to0$, the range of $F$ is relatively compact.
\end{lemma}

\begin{proof}
For each $j$ let $P_j:C(K)\to\ell_\infty^{n_j}$ be the
contractive evaluation map at $n_j$ distinct points of
$K\setminus\{t_*\}$. Choose smooth functions
$\theta_j:\R\to[0,1]$, with pairwise disjoint compact supports
contained in $(0,\infty)$ and tending to infinity, so that there are
$c_j>0$ satisfying
\begin{equation}\label{eq:glue-cutoff}
 \theta_j(c_j)=1,\qquad
 \|\theta'_j\|_\infty\|H_j\|_\infty\le1-b_j^{-1}.
\end{equation}
Dilate a fixed cutoff to obtain the derivative bound, then translate it
beyond the preceding supports. Define $F:C(K)\to X$ by
\begin{equation}\label{eq:glue-definition}
 F(f)=\sum_j\theta_j(f(t_*))\,H_j(P_j f/b_j).
\end{equation}
The sum is locally finite, since $f\mapsto f(t_*)$ is continuous
and the supports escape to infinity. In particular, $F$ is Fr\'echet
$C^\infty$. At each $f$ at most one summand contributes to its
derivative. The product rule, the contraction property of $P_j$ and
$H_j$, and \eqref{eq:glue-cutoff} give
\begin{equation}\label{eq:glue-derivative}
 \|DF(f)\|\le
 \|\theta'_j\|_\infty\|H_j\|_\infty
 +\frac{\theta_j(f(t_*))}{b_j}\le1.
\end{equation}
This remains true at points where all cutoffs vanish. Integrating
\eqref{eq:glue-derivative} along the segment between two functions
in $C(K)$ proves $\Lip(F)\le1$. By normality of $K$, choose pairwise disjointly supported continuous
functions taking the value one at the specified evaluation points,
including $t_*$. The function
$J_j(z)=c_j\psi_*+\sum_{i=1}^{n_j}z_i\psi_i$
then defines an affine isometry from $\ell_\infty^{n_j}$ to $C(K)$
for which $J_j(z)(t_*)=c_j$ and $P_jJ_j(z)=z$.
Equation \eqref{eq:glue-definition} consequently gives
$F(J_j(b_jz))=H_j(z)$, as required. For each fixed $j$, the multiples
$\alpha H_j(z)$, $0\le\alpha\le1$, have compact closure when $H_j$
has relatively compact range. If $\|H_j\|_\infty\to0$, their union is
totally bounded, which proves the final assertion.
\end{proof}

\begin{proof}[Proof of the second assertion of Theorem~\ref{thm:ck}]
Farmer and Johnson proved that
$\pL(A)=\pi_p(A)$ for every bounded linear operator
$A$ \cite[Theorem~2]{FJ}. By the Dvoretzky--Rogers theorem,
$\pi_j(\Id_Y)=\infty$ for every integer $j\ge1$ when $Y$ is
infinite dimensional \cite{DJT}. In
particular, some finite-dimensional subspace $E_j\subset Y$ has
$\pi_j(\Id_{E_j})>2j$. Take the span of a finite family realizing
this strict inequality in $Y$ and use Hahn--Banach on its weak norm.
Theorem~\ref{thm:testing}, applied to $\Id_{E_j}$, gives a contraction
$A_j:\ell_\infty^{n_j}\to E_j$ whose Lipschitz $j$-summing quotient
on a finite test exceeds $2j$.

Clip each source coordinate to an interval containing all its values
on this test. The resulting map is a bounded contraction with the
same test values. Convolution in $\R^{n_j}$ with a nonnegative $C^\infty$
probability kernel supported in a sufficiently small ball preserves
the Lipschitz bound. The resulting maps converge uniformly to the
original bounded contraction. If their uniform difference is at
most $\varepsilon$, Minkowski's inequality bounds the change in
the strong $j$-norm of the selected test by
$2\varepsilon(\sum_i w_i)^{1/j}$; the scalar denominator depends
only on the source points and remains unchanged. Choosing
$\varepsilon$ smaller than the positive margin in the selected test
preserves its strict quotient bound.
Finally, replacing $H(z)$ by $rH(z/r)$ multiplies both the numerator
and scalar supremum on the correspondingly dilated test by $r$.
Consequently, for each integer $j\ge1$, we obtain a bounded,
Fr\'echet $C^\infty$ contraction $H_j$ such that
\begin{equation}\label{eq:large-finite-norm}
 \|H_j\|_\infty\le2^{-j},\qquad
 \pi_j^{\mathrm L}(H_j)>2j.
\end{equation}

Use Lemma~\ref{lem:glue} with $b_j=2$. In view of
\eqref{eq:large-finite-norm}, its contraction $U$ has relatively
compact range. Pulling back the finite test along the affine isometry
$J_j$ and dilating its domain by two shows that
\begin{equation}\label{eq:unbounded-summing}
 \pi_j^{\mathrm L}(U)\ge\frac12\pi_j^{\mathrm L}(H_j)>j.
\end{equation} The inclusion
$\pi_s^{\mathrm L}(U)\le\pL(U)$ for $p\le s$ follows
from the Lipschitz Pietsch domination theorem \cite[Theorem~1]{FJ},
by monotonicity of $L_p$ norms for a probability measure. Taking
integers $j\ge p$ proves the assertion for every finite $p$.
\end{proof}

For a $d$-dimensional real Hilbert space $H_d$, the
$C(K)$ testing formula and the ideal inequality give
\begin{equation}\label{eq:hilbert-growth}
 \sup_{\substack{F:C(K)\to H_d\\ \Lip(F)\le1}}
             \pi_2^{\mathrm L}(F)=\pi_2(\Id_{H_d})=\sqrt d.
\end{equation}
Indeed, for any finite family of vectors, summing their squared
coordinates in an orthonormal basis bounds the strong squared
$2$-norm by $d$ times the weak squared $2$-norm. The basis itself gives
equality. More generally, if $K$ is finite or if $Y$ is finite
dimensional, every Lipschitz map $C(K)\to Y$ is Lipschitz $p$-summing
for each finite $p$. To see this, use the identity on the finite-dimensional
source or target and the ideal inequality for the Lipschitz summing
norm. Its linear summing norm is finite by a decomposition into
rank-one operators.

\begin{corollary}\label{cor:simultaneous}
Let $K$ be an infinite compact Hausdorff space. For each separable
Banach space $X$ there is a single $C^\infty$
contraction $F:C(K)\to X$ such that $\pL(T\circ F)=\pL(T)$ for every
metric space $Y$, every Lipschitz map $T:X\to Y$, and every
$1\le p<\infty$. For any Banach space $X$, a smooth attaining
contraction exists for each fixed $T$ and $p$.
\end{corollary}

\begin{proof}
For $X=\{0\}$ take $F=0$. Otherwise select increasing finite sets
$S_j\subset X$, each containing $0$ and a nonzero point, whose
union is dense in the separable space $X$. In the proof of Theorem~\ref{thm:testing} use the
label set $S_j$ and parameters $\eta=\delta=1/j$. The resulting maps
$A_j:\ell_\infty^{n_j}\to X$ are bounded contractions with range in
the convex hull of $S_j$. The construction of $A_j$ uses only $S_j$
and the two parameters. Convolve to obtain smooth contractions $H_j$
with $\|H_j-A_j\|_\infty\le1/j$. Apply Lemma~\ref{lem:glue} with
$b_j=1+1/j$.

Fix $p,T$ and a finite test in the dense set. For all sufficiently
large $j$ the test lies in $S_j$. Write $W=(\sum_iw_i)^{1/p}$ and
let $A,B$ be its target numerator and scalar supremum. The lifted test
from Theorem~\ref{thm:testing} on the $j$th affine copy has target
numerator at least $A-2\Lip(T)W/j$. Its scalar supremum is at most
$b_j(1+2/j)(B+2C_pW/j)$, by the fiber estimate
\eqref{eq:fiber-oscillation} and the matching inequalities
\eqref{eq:fiber-means}. If $\pL(T\circ F)$ is finite, letting $j\to\infty$ gives
$A\le\pL(T\circ F)B$. Approximating arbitrary endpoints to within
$\varepsilon$ by points in the dense set changes $A$ by at most
$2\Lip(T)\varepsilon W$ and $B$ by at most $2\varepsilon W$. This
proves the reverse inequality. The forward inequality holds by
contractivity, including the case of an infinite norm.

For fixed $T,p$ in an arbitrary Banach space $X$, choose a
sequence of finite tests whose quotients approach $\pL(T)$, with
quotients tending to infinity if the norm is infinite. Their
endpoints belong to a separable closed linear subspace $E\subset X$.
By the McShane extension theorem for real Lipschitz functions,
$\pL(T|_E)=\pL(T)$. Apply the separable case on $E$
and compose the resulting contraction with the inclusion $E\hookrightarrow X$.
\end{proof}

\section{Lattice metrics and Hamming cube lifts}\label{lone:sec:lifts}

A weighted Hamming cube is $\{0,1\}^n$ with distance
$d_w(u,v)=\sum_{j=1}^n w_j|u_j-v_j|$, where $w_j>0$. The fibers in Lemma~\ref{lone:lem:lift} carry constant vector
values and admit measure-preserving matchings.

\begin{lemma}\label{lone:lem:lift}
Let $1\le q\le\infty$, and let $S\subset\Z^d$ be finite, with the
distance inherited from $\ell_q^d$. There exist a weighted Hamming
cube $W$, a contraction $b:W\to\ell_q^d$, and probabilities $\mu_x$ on
$W$, $x\in S$, with the following properties.
\begin{itemize}
\item $b(u)=x$ for every $u\in\supp\mu_x$.
\item For each $x,y\in S$ there is a pair of random vertices $V,V'$
with respective laws $\mu_x,\mu_y$ and
$d_W(V,V')=\norm{x-y}_q$ almost surely.
\end{itemize}
\end{lemma}

\begin{proof}
An integer translation places $S$ in $[0,L]^d$ for some integer
$L\ge2$. Write $C=\Z/(2L)\Z$ with cyclic distance
$d_C(a,b)=\min_{k\in\Z}|a-b+2Lk|$. Embed $C$ in $Q_L$ by starting at
the zero vertex and flipping bits $1,\ldots,L$ in order, twice in
succession. The walk closes after $2L$ steps. A shortest arc has at
most $L$ steps, hence flips distinct bits, so its endpoints have
Hamming distance equal to its length. Denote this isometric embedding
by $\iota$.

Choose nonzero integer columns $v_1,\ldots,v_s$ containing the
standard basis of $\Z^d$ and, up to sign, every nonzero difference
$x-y$ for $x,y\in S$. Put $w_k=\norm{v_k}_q$. The homomorphism
\[
 P:C^s\longrightarrow C^d,\qquad P(u)=\sum_{k=1}^s u_kv_k
\]
is onto because its columns include the standard basis. We view
$C^s$ inside $W=(Q_L)^s$ through the product of the embeddings
$\iota$, and give $W$ the distance
\[
 d_W(u,u')=\sum_{k=1}^s w_k\,\dH(u_k,u'_k).
\]
Take $\mu_x$ to be uniform on $P^{-1}(x)$. If $y-x=\pm v_k$,
translation by $\pm e_k$ is a bijection of $P^{-1}(x)$ onto
$P^{-1}(y)$. It changes only the $k$th cycle coordinate and moves it
one edge, so its displacement in $W$ is $w_k=\norm{x-y}_q$.
For $x=y$ use the identity. This is the coset matching underlying
the transportation identity of Khot and Naor \cite[Lemma~3.6]{KN}.

On the embedded copy of $C^s$ define
\[
 f_j(u)=d_C\big((Pu)_j,0\big),\qquad 1\le j\le d .
\]
Since $0\le x_j\le L$, we have $f_j=x_j$ on $P^{-1}(x)$.
Multiplication by an integer $a$ on a cycle is $|a|$-Lipschitz.
By the triangle inequality on $C$ it follows that
\[
 |f_j(u)-f_j(u')|\le \rho_j(u,u'),\qquad
 \rho_j(u,u')=\sum_{k=1}^s |(v_k)_j|\,\dH(u_k,u'_k).
\]
The same formula defines a pseudometric $\rho_j$ on all of $W$. Set
\[
 b_j(u)=\min_{a\in C^s}\big(f_j(a)+\rho_j(u,a)\big),\qquad u\in W.
\]
The preceding inequality gives $b_j=f_j$ on $C^s$. The triangle
inequality gives $|b_j(u)-b_j(u')|\le\rho_j(u,u')$ on $W$, also when
$\rho_j(u,u')=0$. For $b=(b_1,\ldots,b_d)$ we thus have
\[
 \norm{b(u)-b(u')}_q
 \le \left\|\sum_{k=1}^s \dH(u_k,u'_k)|v_k|\right\|_q
 \le \sum_{k=1}^s \dH(u_k,u'_k)\norm{v_k}_q=d_W(u,u').
\]
Here $|v_k|$ is taken coordinatewise. The first inequality uses
monotonicity of the $\ell_q$ norm and the second is its triangle
inequality, both valid at $q=\infty$. Reversing the initial translation
in the values of $b$ gives the assertion for the original $S$.
\end{proof}

\section{Sampling on metric fibers}\label{lone:sec:sampling}

For a finite metric space $W$ and $m\ge1$, equip $W^m$ with the
averaged sum metric
\[
 d_m(u,v)=\frac1m\sum_{k=1}^m d_W(u_k,v_k).
\]
Changing one coordinate changes a scalar contraction by at most
$\diam(W)/m$, independently of the probability law.

\begin{lemma}\label{lone:lem:sample}
Let $S,W$ be finite metric spaces and $\mu_x$ a probability on $W$
for each $x\in S$. Suppose that for every $x,y$ there is a coupling
$(V,V')$ of $\mu_x,\mu_y$ satisfying $\E d_W(V,V')\le d_S(x,y)$.
Let $Y$ be a Banach space, and let $b:W\to Y$ and $f:S\to Y$
satisfy $b=f(x)$ on $\supp\mu_x$.
Put $b_m(u)=m^{-1}\sum_{k=1}^m b(u_k)$. Then
\[
 \pL(f)\le\liminf_{m\to\infty}\pL(b_m: (W^m,d_m)\to Y)
 \qquad(1\le p<\infty).
\]
\end{lemma}

\begin{proof}
Fix $g\in\Lip_1(W^m,d_m)$ and set
$G(x)=\int g\,d\mu_x^{\otimes m}$. Take $m$ independent copies of a
coupling as in the hypothesis. Their first coordinates form a tuple
$Z$ with law $\mu_x^{\otimes m}$ and their second coordinates form
$Z'$ with law $\mu_y^{\otimes m}$. Consequently
\begin{equation}\label{lone:eq:mean-lipschitz}
 |G(x)-G(y)|\le\E|g(Z)-g(Z')|
 \le\E d_m(Z,Z')\le d_S(x,y).
\end{equation}
By \eqref{lone:eq:mean-lipschitz}, $G$ is an admissible scalar test on $S$.

Write $D=\diam W$. Under $\mu_x^{\otimes m}$, consider the martingale
$M_k=\E[g(Z)\mid Z_1,\ldots,Z_k]$, $0\le k\le m$. With the first
$k-1$ coordinates fixed, the conditional mean as a function of
$Z_k$ has range of length at most $D/m$, by integrating over the
same remaining coordinates. A real random variable supported on an interval of length $a$
has variance at most $a^2/4$. Conditional on the first $k-1$
coordinates, $M_k$ lies in an interval of length $D/m$, and
$M_{k-1}$ is its conditional mean. Orthogonality gives
\begin{equation}\label{lone:eq:variance}
 \Var_{\mu_x^{\otimes m}}g
 =\sum_{k=1}^m\E|M_k-M_{k-1}|^2\le\frac{D^2}{4m}.
\end{equation}

The same bounded-differences argument gives a uniform exponential
tail. If $Z$ has mean zero and is supported on an interval of length
$a$, then
$\log\E e^{tZ}\le t^2a^2/8$ for all $t\in\R$.
Indeed, the second derivative of the logarithmic moment-generating
function is a variance under an exponentially tilted law, bounded by
$a^2/4$; its value and first derivative vanish at zero.
Conditioning successively on $Z_1,\ldots,Z_{k-1}$ therefore gives
\[
 \E\exp\bigl(t(g-G(x))\bigr)
 \le\exp\bigl(t^2D^2/(8m)\bigr).
\]
Markov's inequality, applied to both signs and optimized in $t$, yields
\begin{equation}\label{lone:eq:concentration}
 \mu_x^{\otimes m}\bigl(\{|g-G(x)|>u\}\bigr)
 \le2\exp\bigl(-2mu^2/D^2\bigr)\qquad(u>0).
\end{equation}
For $1\le p\le2$, \eqref{lone:eq:variance} and H\"older's inequality
bound the $L_p$ deviation by $D/(2\sqrt m)$. For $p>2$,
integration of \eqref{lone:eq:concentration} gives
\[
 \int |g-G(x)|^p\,d\mu_x^{\otimes m}
 \le p\int_0^\infty 2u^{p-1}e^{-2mu^2/D^2}\,du
 =2^{1-p/2}\Gamma(1+p/2)\frac{D^p}{m^{p/2}}.
\]
Thus, uniformly in the scalar contraction $g$ and the label $x$,
\begin{equation}\label{lone:eq:deviation}
 \|g-G(x)\|_{L_p(\mu_x^{\otimes m})}\le\delta_{m,p}:=
 \begin{cases}
 D/(2\sqrt m),&1\le p\le2,\\[2pt]
 \bigl(2^{1-p/2}\Gamma(1+p/2)\bigr)^{1/p}D/\sqrt m,&p>2.
 \end{cases}
\end{equation}
For each fixed $p$, $\delta_{m,p}\to0$ as $m\to\infty$,
uniformly in $x$ and $g$.

For a finite test $(x_i,y_i,a_i)_i$ on $S$, take independent tuples
$U_i\sim\mu_{x_i}^{\otimes m}$, $V_i\sim\mu_{y_i}^{\otimes m}$.
Assign the pair $(u,v)$ the weight
$a_i\mu_{x_i}^{\otimes m}(u)\mu_{y_i}^{\otimes m}(v)$ in the $i$th
test. The vector increment is exactly $f(x_i)-f(y_i)$ on every
pair of positive weight. For the scalar increments, add and subtract
$G(x_i)$ and $G(y_i)$ and apply Minkowski's inequality on the weighted
disjoint union of the tests. This gives
\begin{equation}\label{lone:eq:sampling-comparison}
 \left(\sum_i a_i\E|g(U_i)-g(V_i)|^p\right)^{1/p}
 \le \cD_{p,S}\big((x_i,y_i,a_i)_i\big)
     +2\delta_{m,p}\left(\sum_i a_i\right)^{1/p}.
\end{equation}
In particular, \eqref{lone:eq:sampling-comparison} gives the finite-$m$
estimate
\begin{align*}
 &\left(\sum_i a_i\|f(x_i)-f(y_i)\|^p\right)^{1/p}\le\pL(b_m)\biggl[
  \cD_{p,S}\big((x_i,y_i,a_i)_i\big)
     +2\delta_{m,p}\Big(\sum_i a_i\Big)^{1/p}\biggr].
\end{align*}
Letting $m\to\infty$ along a subsequence attaining the lower limit
proves the assertion when that limit is finite. The case of infinite
lower limit is immediate.
\end{proof}

Parallel coordinate segments in an $\ell_{p'}$ lattice box satisfy
a dimension-free scalar estimate, up to the lattice approximation
error.

\begin{lemma}\label{lone:lem:box}
Let $1\le p<\infty$, $q=p'$, and take integers $1\le h<L$. Put
$\Lambda=\{0,\ldots,L\}^d$ with its $\ell_q$ metric. For each
$i=1,\ldots,d$ take the pairs $(a,a+he_i)$ with
\[
 a\in\{0,\ldots,L-1\}^d,\qquad a_i<L-h,
\]
each with weight $L^{-d}$. The total weight in each direction is
$\lambda=1-h/L$, and the scalar supremum for this test is at most
\[
 h+2d^{1/q}(d\lambda)^{1/p}.
\]
At $p=1$ we use $q=\infty$ and $d^{1/q}=1$.
\end{lemma}

\begin{proof}
Each admissible $g$ on $\Lambda$ extends to a scalar contraction
on $\ell_q^d$. Use the same symbol for the extension. Let
$V_i=\{x\in[0,L)^d:x_i<L-h\}$ and take coordinatewise integer parts.
For $x\in V_i$, $\lfloor x+he_i\rfloor=\lfloor x\rfloor+he_i$ and
$\|x-\lfloor x\rfloor\|_q\le d^{1/q}$. Replacing the integer
increment by $g(x+he_i)-g(x)$ therefore costs at most $2d^{1/q}$
pointwise. Integration with measure $dx/L^d$ counts each integer
starting point with weight $L^{-d}$. The $L_p$ cost over all $V_i$
is at most $2d^{1/q}(d\lambda)^{1/p}$.

At every point of differentiability, the differential of $g$ has
norm at most one on
$\ell_q^d$, hence $\sum_i|\partial_i g|^p\le1$. Rademacher's theorem
gives differentiability almost everywhere. Absolute continuity on
coordinate lines and H\"older's inequality give, for almost every
$x\in V_i$,
\[
 |g(x+he_i)-g(x)|^p
 \le h^{p-1}\int_0^h |\partial_i g(x+te_i)|^p\,dt .
\]
For $0\le t\le h$, the translate $V_i+te_i$ is contained in $[0,L]^d$.
Fubini's theorem and the gradient bound give
\begin{equation}\label{lone:eq:box-continuous}
 \sum_{i=1}^d\int_{V_i}|g(x+he_i)-g(x)|^p\,\frac{dx}{L^d}
 \le h^p\int_{[0,L]^d}\sum_{i=1}^d|\partial_i g(x)|^p\,\frac{dx}{L^d}
 \le h^p .
\end{equation}
Taking $p$th roots and restoring the lattice error proves the lemma.
\end{proof}

\section{Recovery of the absolutely summing norm}\label{lone:sec:recovery}

Multilinear interpolation preserves the Lipschitz bound on a
finite weighted Hamming cube.

For $A:\{0,1\}^n\to X$, define its multilinear interpolant on
$[0,1]^n$ by
\[
 \widetilde A(t)=\sum_{\varepsilon\in\{0,1\}^n}
 A(\varepsilon)\prod_{\varepsilon_j=1}t_j
                    \prod_{\varepsilon_j=0}(1-t_j).
\]
If parallel edge increments in direction $j$ have norm at most
$w_j$, then $\|\partial_j\widetilde A(t)\|\le w_j$. Indeed, the
derivative is a convex combination of those increments, with
coefficients given by the same product over the other $n-1$
coordinates. Varying one coordinate at a time gives
\[
 \|\widetilde A(t)-\widetilde A(t')\|
 \le \sum_{j=1}^n w_j|t_j-t'_j|.
\]
Composing with coordinatewise clipping to $[0,1]$ gives an extension
to the entire weighted $\ell_1^n$ space with the same bound. Its
range lies in the convex hull of the finitely many vertex values.

\begin{proof}[Proof of Theorem~\ref{lone:thm:recovery}]
Let $K$ be the supremum over ordinary Hamming cubes in the theorem.
If $F:M\to X$ is a contraction on any metric space, apply the linear
summing inequality to the vectors
$a_i^{1/p}(F(s_i)-F(t_i))$. The function $x^*\circ F$ belongs to
$\Lip_1(M)$ whenever $x^*\in B_{X^*}$, so
$\pL(U\circ F)\le\pi_p(U)$. We prove $\pi_p(U)\le K$ and then construct a contraction on
$\ell_1$ attaining this value. For the first assertion we may assume
$K<\infty$.

Let $A:Q_n\to X$ be a contraction when the $j$th bit has weight
$w_j>0$. Choose $\tau>0$, let $k_j=\lceil w_j/\tau\rceil$, and
write $w'_j=\tau k_j$. Replace bit $j$ by a block of $k_j$ equal
bits. On the larger cube, of dimension $k_1+\cdots+k_n$, define
a map by taking the arithmetic mean of each block and applying
$\widetilde A$ to the resulting $n$ numbers. Flipping a bit in
block $j$ changes its value by at most $w_j/k_j\le\tau$.
It agrees with $A$ on the repeated-bit vertices. On the larger cube with metric $\tau\dH$ the resulting map $B$
is a contraction. Since $U$ is linear, $B/\tau$ is a contraction on
the ordinary Hamming cube. On every finite weighted test,
changing the metric from $\tau\dH$ to $\dH$ divides the
scalar denominator by $\tau$; replacing $U\circ B$ with
$U\circ (B/\tau)$ divides the vector numerator by the
same factor. Therefore
\[
 \pL(U\circ B:(Q_M,\tau\dH)\to Y)
 =\pL(U\circ(B/\tau):(Q_M,\dH)\to Y)
 \le K.
\]
Here $M=k_1+\cdots+k_n$. On the repeated-bit set the ambient metric induces
$d_{w'}$. Put $c=\max_j w'_j/w_j\ge1$. Then
$d_w\le d_{w'}\le c d_w$. For any weighted test of pairs on
$Q_n$, its vector numerator for $U\circ A$ does not depend on
which of these two metrics is used. By restriction and the McShane
extension theorem, the corresponding scalar denominator for
$d_{w'}$ is the same whether computed on the repeated-bit
subset or on the larger cube. Moreover, every scalar contraction
for $d_{w'}$ is $c$-Lipschitz for $d_w$; thus
$\mathcal D_{p,d_{w'}}\le c\mathcal D_{p,d_w}$ for each fixed
test. The bound $K$ on the larger cube therefore implies
\begin{equation}\label{lone:eq:weighted-bound}
 \pL(U\circ A:(Q_n,d_w)\to Y)\le K\max_j\frac{w'_j}{w_j}.
\end{equation}
Letting $\tau\downarrow0$ in \eqref{lone:eq:weighted-bound} gives
$\pL(U\circ A)\le K$ on every finite weighted Hamming cube.

Now fix $u_1,\ldots,u_d\in X$, take $q=p'$, and let
\[
 w_p(u_1,\ldots,u_d)
 =\sup_{x^*\in B_{X^*}}
       \left(\sum_{i=1}^d|x^*(u_i)|^p\right)^{1/p}
\]
be the weak $p$-norm of the family. Define
$B:\ell_q^d\to X$ by $Be_i=u_i$. Duality gives
\begin{equation}\label{lone:eq:B-norm}
 \norm B
 =\sup_{x^*\in B_{X^*}}\sup_{\norm c_q\le1}
       \left|\sum_i c_i x^*(u_i)\right|
 =w_p(u_1,\ldots,u_d).
\end{equation}
For each finite $S\subset\Z^d$ take $W,b,\mu_x$ from
Lemma~\ref{lone:lem:lift}. The lemma supplies couplings with
$\mathbb E d_W(V,V')=\|x-y\|_q$ for all $x,y\in S$;
moreover, $UBb$ equals $UBx$ on $\operatorname{supp}\mu_x$.
Thus all hypotheses of Lemma~\ref{lone:lem:sample} hold for the
map $x\mapsto UBx$ on $S$ and the map $UBb$ on $W$. For every $m$ the
averaged map on $W^m$ is $U$ composed with
\[
 (z_1,\ldots,z_m)\longmapsto\frac1m\sum_{k=1}^m B b(z_k).
\]
Its inner map has Lipschitz constant at most $\norm B$ for $d_m$.
This metric is a weighted Hamming metric, with the original bit
weights divided by $m$ in each copy. The weighted-cube estimate \eqref{lone:eq:weighted-bound},
together with homogeneity, gives
\begin{equation}\label{lone:eq:lattice-recovery}
 \pL(UB|_S)\le K\norm B .
\end{equation}

Apply this inequality to the lattice test of Lemma~\ref{lone:lem:box}. On
every pair in direction $i$ the vector increment has norm
$h\|Uu_i\|$. With $\lambda=1-h/L$, we obtain
\begin{equation}\label{lone:eq:lattice-estimate}
 h\lambda^{1/p}\left(\sum_{i=1}^d\|Uu_i\|^p\right)^{1/p}
 \le K\norm B\,\big(h+2d^{1/q}(d\lambda)^{1/p}\big).
\end{equation}
In \eqref{lone:eq:lattice-estimate}, set $L=h^2$, divide by $h$,
and let $h\to\infty$ through the integers with $d$ fixed. By
\eqref{lone:eq:B-norm}, this gives
$(\sum_i\|Uu_i\|^p)^{1/p}\le K w_p(u_1,\ldots,u_d)$ for every
finite family, hence $\pi_p(U)\le K$.

Finally, choose cube contractions $A_j:Q_{n_j}\to X$ with
$\pL(U\circ A_j)$ tending to $K$, also when $K=\infty$.
Subtract $A_j(0)$ and apply the clipped multilinear extension to
obtain contractions $F_j:\ell_1^{n_j}\to X$ with $F_j(0)=0$.
Identify $\ell_1$ with the $\ell_1$-sum of these coordinate blocks
and set
\[
 F((z_j)_{j\ge1})=\sum_{j\ge1}F_j(z_j).
\]
The series converges absolutely since $F_j(0)=0$ and
$\sum_j\|F_j(z_j)\|\le\sum_j\|z_j\|_1$. For $z,z'$,
\[
 \|F(z)-F(z')\|\le\sum_j\|F_j(z_j)-F_j(z'_j)\|
 \le\sum_j\|z_j-z'_j\|_1=\|z-z'\|_1.
\]
Its restriction to each coordinate block recovers $F_j$.
Moreover, every real scalar contraction on an embedded cube
extends to $\ell_1$ with the same Lipschitz constant, so the
summability norm of the restricted map is preserved. Consequently
$\pL(U\circ F)\ge\sup_j\pL(U\circ A_j)=K$.
\end{proof}

\section{Smooth counterexamples on \texorpdfstring{$L_1$}{L1}}
\label{lone:sec:consequences}

For $1\le p\le s<\infty$, the monotonicity inequality
$\pL(T)\ge\pi_s^{\mathrm L}(T)$
\cite[Section~2]{FJ} allows a single smooth map to fail Lipschitz
$p$-summability for every finite $p$. For $s>p$, this also
follows directly from the definition. On the finite measure space
with masses $a_i$, take nonnegative multipliers $w$ of norm at most
one in $L_{s/(s-p)}$. Apply the $p$-summing inequality with weights
$a_i w_i$. H\"older bounds the scalar side by the scalar $s$-norm
for the original weights. Taking the supremum over $w$ gives the
$s$-norm on the vector side by duality.

\begin{proof}[Proof of Theorem~\ref{lone:thm:failure}]
The Dvoretzky--Rogers theorem gives $\pi_j(\Id_Y)=\infty$ for every
integer $j\ge1$ when $Y$ is infinite dimensional \cite{DJT}. By
Theorem~\ref{lone:thm:recovery} there is a cube contraction
$A_j:Q_{n_j}\to Y$ and a finite test on which its Lipschitz
$j$-summing inequality fails with constant $j$.

The clipped interpolant $G_j:\ell_1^{n_j}\to Y$ constructed in
Section~\ref{lone:sec:recovery} is bounded and contractive, with range in
the span of the values of $A_j$. Let $\rho_\varepsilon$ be a
nonnegative smooth mollifier of mass one supported in an
$\ell_1^{n_j}$ ball of radius $\varepsilon$. Convolution preserves
the Lipschitz bound and the finite-dimensional range space.
Moreover $\|G_j*\rho_\varepsilon-G_j\|_\infty\le\varepsilon$.
On the chosen test the vector $j$-norm changes by at most
$2\varepsilon(\sum_i a_i)^{1/j}$, while the scalar supremum is
unchanged. Thus for sufficiently small $\varepsilon$ we obtain a
bounded $C^\infty$ contraction $H_j$ with
$\pi_j^{\mathrm L}(H_j)>j$.

Fix $\eta\in C_c^\infty((-1/3,1/3))$ with $0\le\eta\le1$ and
$\eta(0)=1$. The simultaneous dilation
$z\mapsto r H_j(z/r)$ preserves both the Lipschitz and the summing
norms. Indeed, rescaling the endpoints by $r$ multiplies both
sides of a finite test by $r$, and the inverse rescaling gives the
reverse comparison. Choosing $r$ separately for each block, we
may impose
\[
 \|H_j\|_\infty\le\frac{2^{-j}}{1+\|\eta'\|_\infty}.
\]
Identify $\ell_1$ with
$\R\oplus_1(\bigoplus_{j\ge1}\ell_1^{n_j})_{\ell_1}$. Define
\[
 T(t,z)=\sum_{j\ge1}\eta(t-j)H_j(z_j).
\]
The supports in $t$ are disjoint and locally finite. Each term
depends smoothly on finitely many coordinates, so $T$ is Fr\'echet
$C^\infty$. Where the $j$th term is active its derivative on
$(s,v)$ is
\[
 s\,\eta'(t-j)H_j(z_j)+\eta(t-j)DH_j(z_j)v_j .
\]
Its norm is at most $|s|+\|v_j\|_1$. Elsewhere the derivative is
zero. Integration along line segments shows that $T$ is a contraction.
The values of each term lie in a bounded subset of a finite-dimensional
space. All terms of index larger than $J$ take values of norm at most
$2^{-J}$. Since only one term can be nonzero at a given point,
the range is totally bounded, hence relatively compact.

On the isometric copy of $\ell_1^{n_j}$ with $t=j$ and all other
blocks zero, $T$ coincides with $H_j$. Thus for every fixed finite
$p\ge1$ and every integer $j\ge p$,
\begin{equation}\label{lone:eq:simultaneous-failure}
 \pL(T)\ge\pi_j^{\mathrm L}(T)
 \ge\pi_j^{\mathrm L}(H_j)>j.
\end{equation}
Since $j$ is arbitrary, \eqref{lone:eq:simultaneous-failure} gives
$\pL(T)=\infty$ for every $1\le p<\infty$.
\end{proof}

\begin{corollary}\label{lone:cor:classification}
Let $X=L_1(\mu)$, let $Y$ be a Banach space, and fix $1\le p<\infty$.
Every Lipschitz map $X\to Y$ is Lipschitz $p$-summing if and only if
at least one of the spaces is finite dimensional. If both are
infinite dimensional, there is a $C^\infty$ contraction with
relatively compact range for which every finite Lipschitz summing
norm is infinite.
\end{corollary}

\begin{proof}
If $L_1(\mu)$ is infinite dimensional, there are countably many
disjoint sets $E_j$ of finite positive measure. Indeed, infinitely
many distinct finite-measure
atoms give such a family after removal of null intersections. If
there are only finitely many such atoms up to null sets, either they
support all integrable functions, giving a finite-dimensional space,
or a finite positive-measure set remains outside their union. Such a
set is nonatomic and can be split successively, setting one positive
piece aside at each step. Here it suffices to consider sets of finite
measure, since $\{|f|>t\}$ has finite measure for every
$f\in L_1(\mu)$ and $t>0$.

With $u_j=\mathbf1_{E_j}/\mu(E_j)$, the map
$J(c)=\sum_j c_j u_j$ is an isometry from $\ell_1$ into $L_1(\mu)$.
Its left inverse $R f=(\int_{E_j}f\,d\mu)_j$ has norm at most one,
by disjointness. Precomposition by $R$ preserves the summing norm of
every map on $\ell_1$, as precomposition by $J$ recovers that map.
Applied to the map of Theorem~\ref{lone:thm:failure}, this
precomposition also preserves smoothness and the range.

The identity of a finite-dimensional space has finite $\pi_p$, by
decomposition into finitely many rank-one operators. If $X$ is
finite dimensional, bounding each vector increment of $T$ by
$\Lip(T)$ times the source distance gives
$\pL(T)\le\Lip(T)\pL(\Id_X)\le\Lip(T)\pi_p(\Id_X)$.
If $Y$ is finite dimensional, applying its linear summing inequality
to the increments of $T$ gives $\pL(T)\le\pi_p(\Id_Y)\Lip(T)$.
\end{proof}

\begin{corollary}\label{lone:cor:rank}
If $H$ is a real Hilbert space of dimension $d\ge1$, then for
$1\le p<\infty$,
\[
 \max_{\substack{T:\ell_1\to H\\ \Lip(T)\le1}}\pL(T)
 =\left(
 \frac{\sqrt\pi\,\Gamma((d+p)/2)}
      {\Gamma(d/2)\,\Gamma((p+1)/2)}
 \right)^{1/p}.
\]
The same value holds on every infinite-dimensional $L_1(\mu)$.
\end{corollary}

\begin{proof}
By Theorem~\ref{lone:thm:recovery} it suffices to compute
$\pi_p(\Id_H)$. The classical constant has the following form
\cite[p.~291]{Pietsch}. Let $\sigma$ be uniform probability on
$S^{d-1}$, and put $c^p=\int|\theta_1|^p\,d\sigma(\theta)$.
Rotation invariance gives
$\int|\langle v,\theta\rangle|^p\,d\sigma=c^p\|v\|_2^p$.
Summing over any finite family and bounding the integral by its
supremum gives $\pi_p(\Id_H)\le c^{-1}$.

Conversely, partition the sphere into finitely many sets of diameter
at most $\varepsilon$. Let $a_i$ be their masses and choose a point
$u_i$ in each nonempty cell. Uniformly in $\theta\in S^{d-1}$,
\[
 \left|\sum_i a_i|\langle u_i,\theta\rangle|^p-c^p\right|
 \le p\varepsilon ,
\]
since $t\mapsto |t|^p$ is $p$-Lipschitz on $[-1,1]$.
The family $a_i^{1/p}u_i$ has strong $p$-norm one and weak
$p$-norm at most $(c^p+p\varepsilon)^{1/p}$. Letting
$\varepsilon\downarrow0$ yields $\pi_p(\Id_H)\ge c^{-1}$.

For a standard Gaussian vector $Z$ in $\R^d$, the direction is
uniform and independent of its length, so
$c^p=\E|Z_1|^p/\E\|Z\|_2^p$. Integration in polar coordinates,
with the substitution $s=r^2/2$, gives
\[
 \E|Z_1|^p=
 \frac{2^{p/2}\Gamma((p+1)/2)}{\sqrt\pi},
 \qquad
 \E\|Z\|_2^p=
 \frac{2^{p/2}\Gamma((d+p)/2)}{\Gamma(d/2)} .
\]
Their ratio gives the formula, including $d=1$.
\amend{For an infinite-dimensional $L_1(\mu)$, the isometry $J$ in
Corollary~\ref{lone:cor:classification} has a norm-one left inverse $R$,
so $JR$ projects onto $J(\ell_1)$. Since $RJ=\Id_{\ell_1}$, precomposition with $R$
preserves the attained value. The upper bound follows from
$\pL(T)\le\pi_p(\Id_H)\Lip(T)$.}
\end{proof}

\section{Hilbert-valued maps on the Hamming cube}\label{lone:sec:cube}

Let $K_N$ be the supremum of $\pi_1^{\mathrm L}(F)$ over all
contractions from $Q_N$ into a real Hilbert space. Each map has range in an affine space of dimension at most
$2^N-1$. Translation-invariant pair families have uniformly bounded
summing quotients, although $K_N$ grows with $N$.

\begin{proposition}\label{lone:prop:invariant}
Identify $Q_N$ isometrically with $\{-1,1\}^N$, retaining the
Hamming metric, and let $v^S$ denote
$v$ with the signs in $S\subset\{1,\ldots,N\}$ reversed. For any
Lipschitz $F:Q_N\to H$ and any $a_S\ge0$,
\[
 \sum_S a_S\,\E_v\|F(v)-F(v^S)\|
 \le\sqrt2\,\Lip(F)\sup_{g\in\Lip_1(Q_N)}
       \sum_S a_S\,\E_v|g(v)-g(v^S)| ,
\]
where $v$ is uniform on $Q_N$. The constant $\sqrt2$ is sharp for
every $N\ge2$.
\end{proposition}

\begin{proof}
Expand $F=\sum_A\widehat F(A)\chi_A$ in the Walsh characters
$\chi_A(v)=\prod_{i\in A}v_i$. Parseval's identity applies
coordinatewise in the finite-dimensional span of the range. Since
$\chi_A(v^S)=(-1)^{|A\cap S|}\chi_A(v)$,
\[
 \begin{aligned}
 \E_v\|F(v)-F(v^S)\|^2
 &=4\sum_{|A\cap S|\ {\rm odd}}\|\widehat F(A)\|^2\\
 &\le\sum_{i\in S}\E_v\|F(v)-F(v^{\{i\}})\|^2
 \le |S|\,\Lip(F)^2 .
 \end{aligned}
\]
Every Fourier coefficient on the first line is counted at least
once on the second. Cauchy--Schwarz gives the corresponding
first-moment bound $\sqrt{|S|}\,\Lip(F)$.

The scalar contraction $g(v)=\frac12\sum_i v_i$ has
$\E_v|g(v)-g(v^S)|=\E|\sum_{i\in S}v_i|\ge\sqrt{|S|/2}$.
Put $\alpha_k=\E|\sum_{i=1}^k v_i|$.
Conditioning on the first $k$ signs gives
$\alpha_{k+1}-\alpha_k=\Pr(\sum_{i=1}^k v_i=0)$.
It follows by induction that
\[
 \alpha_{2m-1}=\alpha_{2m}=2m\binom{2m}{m}4^{-m},\qquad m\ge1 .
\]
The successive ratios of $\alpha_{2m}/\sqrt{2m}$ are
$(2m+1)/(2\sqrt{m(m+1)})\ge1$, and the initial value is
$1/\sqrt2$. This proves the scalar bound for even and odd $k$.
Multiplying by $a_S$ and summing proves the inequality.

For sharpness take the two diagonals of $Q_2$ and $F(v)=v/2$
in $\R^2$. Both diagonals have length $\sqrt2$. For any scalar
contraction, orient each diagonal from its greater value to its
lesser value. The sum of the two decrements is unchanged if each
greater endpoint is paired with the lesser endpoint of the other
diagonal. Those two new pairs are edges, so the sum is at most two.
The function $(v_1+v_2)/2$ attains two, giving scalar average one.
In higher dimensions take $S=\{1,2\}$ and use the same example on
each two-dimensional coordinate slice.
\end{proof}

\begin{theorem}\label{lone:thm:rate}
For every integer $r\ge1$, set $N_r=(8r\,3^r)^3$. There is a
contraction $Q_{N_r}\to\ell_2^r$ whose Lipschitz $1$-summing norm
is at least $4\sqrt r/11$. Consequently, for an absolute $c>0$,
\begin{equation}\label{lone:eq:cube-rate}
 c\sqrt{\log(N+1)}\le K_N\le\sqrt{2N}\qquad(N\ge1).
\end{equation}
\end{theorem}

\begin{proof}
Take $L=8r$, $h=L/2$, $s=3^r$ and $D=Ls$. In the cyclic construction
of Lemma~\ref{lone:lem:lift} take as columns the $s$ vectors of
$\{-1,0,1\}^r$, assigning every block weight one, including the
unused zero column. Thus $W=Q_D$ with ordinary Hamming distance.
Use the same coordinate functions and extensions as in that lemma.
On any single edge of $W$, every coordinate changes by at most one,
so $b:W\to\ell_\infty^r$ is a contraction. All labels
$a\in\Lambda=\{0,\ldots,L\}^r$ have uniform fiber measures
$\mu_a$ with $b=a$ on their supports.

For these columns the coupling estimate holds for every two lattice
points. If $a'-a=z$ and $k=\|z\|_\infty$, write
$z=\sum_{\ell=1}^k v^{(\ell)}$, where
\[
 v^{(\ell)}_j=\operatorname{sgn}(z_j)\,
                   \mathbf1_{\{|z_j|\ge\ell\}},\qquad 1\le j\le r .
\]
Each summand is one of the columns. Composing the corresponding $k$
unit translations takes the fiber of $a$ bijectively onto that of
$a'$ and displaces every representative by at most $k$. In particular
the means of scalar contractions under the product measures are
$1$-Lipschitz in the $\ell_\infty$ metric on $\Lambda$.

Take $m=D^2$ and define on $(W^m,d_m)$ the map
\[
 T(u_1,\ldots,u_m)=\frac1{m\sqrt r}\sum_{j=1}^m b(u_j)\in\ell_2^r .
\]
It is a contraction since $\|z\|_2\le\sqrt r\,\|z\|_\infty$.
On $\supp\mu_a^{\otimes m}$ it equals $a/\sqrt r$. Use the test of
Lemma~\ref{lone:lem:box} with $p=1$, $d=r$, and lift each pair of
labels by independent tuples from the two product measures. There
is total mass $1/2$ in each direction, and every vector increment
has length $h/\sqrt r$. The numerator is therefore $h\sqrt r/2$.

The means of the scalar tests have total increment at most $h+r$ by
Lemma~\ref{lone:lem:box}. A sample differs from its mean in $L_1$ by at
most $D/(2\sqrt m)=1/2$, by the variance estimate in
Lemma~\ref{lone:lem:sample}. Returning from the two means to their
samples costs at most one per pair, or $r/2$ in total. Thus
\[
 \pi_1^{\mathrm L}(T)\ge\frac{h\sqrt r/2}{h+3r/2}
                         =\frac4{11}\sqrt r .
\]
The metric $d_m$ is $1/m$ times ordinary Hamming distance on
$Q_{Dm}$. Hence $mT$ is a contraction on the ordinary cube and has
the same summing norm that $T$ has for $d_m$. Here $Dm=D^3=N_r$.
The sequence $K_N$ is nondecreasing by projection onto coordinate
subcubes. For $N_r\le N<N_{r+1}$ use
$\log(N_{r+1}+1)\le C r$ with an absolute $C$.
The initial dimensions are covered by $K_N\ge1$.

For the upper bound in \eqref{lone:eq:cube-rate}, choose independent random
signs $\varepsilon_i$ and take the scalar contraction
$g_\varepsilon(v)=\frac12\sum_i\varepsilon_i v_i$ in sign
coordinates. On a pair of distance $k$ its expected absolute
increment is $\alpha_k\ge\sqrt{k/2}$, by the calculation in
Proposition~\ref{lone:prop:invariant}. The scalar supremum is at least
the average over these contractions. A Hilbert-valued contraction
has increment at most $k\le\sqrt{2N}\sqrt{k/2}$ on the same pair.
Multiplication by the weights and summation prove the upper bound.
\end{proof}

\section{A summing map without an integral factorization}\label{sec:integral}

Let $Y$ be a pointed Banach space, and write
$\operatorname{Lip}_0(Y)$ for the Banach space of real-valued
Lipschitz functions vanishing at zero, equipped with the Lipschitz
norm. The canonical map
$\delta_Y:Y\to\operatorname{Lip}_0(Y)^*$ is defined by
$\delta_Y(y)(\varphi)=\varphi(y)$. Following \cite{FJ}, a pointed
Lipschitz map $T:M\to Y$ is Lipschitz $q$-integral if there are a
probability space and Lipschitz maps $A$ and $B$ for which
\begin{equation}\label{eq:integral-factorization}
 \delta_Y\circ T=B\circ i_{\infty,q}\circ A,\qquad
 A:M\to L_\infty(\mu),\quad
 B:L_q(\mu)\to\operatorname{Lip}_0(Y)^*.
\end{equation}
Here $i_{\infty,q}:L_\infty(\mu)\to L_q(\mu)$ is the canonical
inclusion. \amend{Ch\'avez-Dom\'inguez and Jim\'enez-Vargas
\cite[Proposition~2.7(b)]{CDJV} proved that a Lipschitz $q$-integral
map into a dual Banach space extends to every metric superspace of its
domain with the same Lipschitz $q$-integral norm.}
For countable domains, a factorization of the form
\eqref{eq:integral-factorization} yields a Lipschitz extension to
the bidual target, as shown in the proof of
Theorem~\ref{thm:integral}.

\amend{The target norm in Lemma~\ref{lem:related-integral-block}
adapts Kalton's method of sections \cite{Kalton}. Mendel and Naor used
this construction to obtain an extension lower bound of order $d^{1/4}$
for maps on nets in Euclidean $d$-space
\cite[\S6.1 and Theorem~1.17]{MN}.}
A net of the Euclidean sphere, separated in $\ell_1^d$, gives
finite examples whose image distances are controlled by two families
of scalar contractions. Grothendieck's inequality then bounds the trace
of the derivative of any extension.

\begin{lemma}\label{lem:related-integral-block}
For every $d\ge16$ there are a finite set $M_d\subset\ell_\infty^{2^d}$,
a finite-dimensional Banach space $Y_d$, and a map $f_d:M_d\to Y_d$
with $\pi_1^{\mathrm L}(f_d)\le4$. Every extension of $f_d$ to
$\ell_\infty^{2^d}$ with values in $Y_d$ has Lipschitz constant at
least $c d^{1/4}$, where $c>0$ is absolute.
\end{lemma}

\begin{proof}
Put $\delta=d^{-1/4}$ and choose a maximal $\delta$-separated set $S$
on the Euclidean unit sphere, using the $\ell_1^d$ distance. By
compactness $S$ is finite and every point of the sphere is within
$\ell_1^d$ distance $\delta$ of $S$. Let $\Omega=\{-1,1\}^d$, equipped with the uniform
probability measure, and define
$R:\R^d\to\R^\Omega$ by $(Rx)_r=\langle r,x\rangle$.
Since $\max_{r\in\Omega}|\langle r,x\rangle|=\|x\|_1$,
$R:\ell_1^d\to\ell_\infty^\Omega$ is a linear isometry.
Set $M_d=R(S)$. On $\R^S$ let $Q:\R^S\to\R^d$ be the linear map
specified by $Qe_s=s$, and equip $\R^S$ with the norm
\begin{equation}\label{eq:finite-target-norm}
 \|u\|_{Y_d}=\max\bigl\{(\delta/2)\|u\|_1,\ \|Qu\|_2\bigr\}.
\end{equation}
Define $f_d:M_d\to Y_d$ by $f_d(Rs)=e_s$.
The separation of $S$ ensures that $f_d$ is a contraction.
For distinct $s,t$ the image distance is
$\max\{\delta,\|s-t\|_2\}$.

Consider the following probability distribution on scalar functions
$g:M_d\to\R$. With probability $1/2$, choose
$g(Rs)=\langle r,s\rangle$, where $r\in\Omega$ is uniform.
With the remaining probability $1/2$, choose
$g(Rs)=(\delta/2)\xi_s$, where the $\xi_s$ are independent
Rademacher signs. Each realization is a contraction: for the first
family this follows from the isometry $R$, and for the second from
the $\delta$-separation of $S$.

\amend{Fix distinct $s,t\in S$.} For the first family, put
$Z=\langle r,s-t\rangle$. Independence
of the signs gives
\[
 \E Z^2=\|s-t\|_2^2,\qquad
 \E Z^4=3\|s-t\|_2^4-2\sum_{k=1}^d(s_k-t_k)^4
 \le3\|s-t\|_2^4.
\]
H\"older's inequality, applied to
$|Z|^2=|Z|^{2/3}|Z|^{4/3}$, implies
\begin{equation}\label{eq:rademacher-first-moment}
 \E|Z|\ge\frac{(\E Z^2)^{3/2}}{(\E Z^4)^{1/2}}
 \ge\frac{\|s-t\|_2}{\sqrt3}.
\end{equation}
\amend{Independence} of $\xi_s,\xi_t$ gives
$\E|\xi_s-\xi_t|=1$, so the second family has expected absolute
increment $\delta/2$. Combining the two choices and using
\eqref{eq:finite-target-norm}, we obtain
\begin{equation}\label{eq:random-contractions}
 \E|g(Rs)-g(Rt)|\ge
 \max\left\{\frac{\|s-t\|_2}{2\sqrt3},\frac\delta4\right\}
 \ge\frac14\|e_s-e_t\|_{Y_d}.
\end{equation}
\amend{On the diagonal $s=t$, the expected increment and the image
distance both vanish.} Multiplying by the test weights and summing proves the asserted bound
on $\pi_1^{\mathrm L}(f_d)$.

Let $F:\ell_\infty^\Omega\to Y_d$ be an extension with Lipschitz
constant $L$. Convolution with a smooth probability density supported
in a sufficiently small ball gives a smooth map $F_\varepsilon$ with
$\Lip(F_\varepsilon)\le L$ and $\|F_\varepsilon-F\|_\infty\le\varepsilon$.
For $x$ on the Euclidean unit sphere choose $s\in S$ with
$\|x-s\|_1\le\delta$. Since $Qf_d(Rs)=s$ and $\|Q\|_{Y_d\to\ell_2^d}\le1$,
\[
 \|QF_\varepsilon(Rx)-x\|_2\le (L+1)\delta+\varepsilon.
\]
At each point put $A=DF_\varepsilon(Rx)$. As a matrix from
$\ell_\infty^\Omega$ to $\ell_1^S$, it has norm at most $2L/\delta$.
Write $A_{sr}$ for the matrix coefficients of
$A:\ell_\infty^\Omega\to\ell_1^S$ in the standard bases.
Since $\|s\|_2=\|r/\sqrt d\|_2=1$, the real Grothendieck inequality \cite{Groth}, with constant
$K_G$, gives
\begin{align}
 |\operatorname{tr}(QAR)|
 &=\left|\sum_{s\in S}\sum_{r\in\Omega}
                  A_{sr}\langle s,r\rangle\right|\notag\\
 &\le K_G\sqrt d\sup_{\substack{|\alpha_s|\le1\\|\beta_r|\le1}}
                \left|\sum_{s,r}A_{sr}\alpha_s\beta_r\right|\notag\\
 &=K_G\sqrt d\,\|A\|_{\ell_\infty^\Omega\to\ell_1^S}.
 \label{eq:groth-trace}
\end{align}
The last equality is the dual characterization of the indicated
operator norm. By \eqref{eq:finite-target-norm} and
$\|DF_\varepsilon\|\le L$, the final expression in
\eqref{eq:groth-trace} is at most $2K_GL\sqrt d/\delta$.
Let $B_2^d$ denote the Euclidean unit ball and
$\sigma_{d-1}$ the surface measure of its boundary. The map
$\Phi_\varepsilon=QF_\varepsilon R:B_2^d\to\R^d$ is smooth, and
for $x\in\partial B_2^d$ we have
$\langle\Phi_\varepsilon(x),x\rangle
\ge1-(L+1)\delta-\varepsilon$.
The divergence theorem and \eqref{eq:groth-trace} therefore give
\begin{align}
 &\bigl(1-(L+1)\delta-\varepsilon\bigr)
   \sigma_{d-1}(\partial B_2^d)\notag\\
 &\qquad\le\int_{B_2^d}\operatorname{div}\Phi_\varepsilon(x)\,dx
 \le \frac{2K_GL\sqrt d}{\delta}\,|B_2^d|.
 \label{eq:divergence}
\end{align}
Since $\sigma_{d-1}(\partial B_2^d)=d|B_2^d|$,
\eqref{eq:divergence} reduces to
\begin{equation}\label{eq:finite-extension-lower}
 1-(L+1)\delta-\varepsilon\le\frac{2K_GL}{\delta\sqrt d}.
\end{equation}
Let $\varepsilon\downarrow0$, use $\delta=d^{-1/4}\le1/2$,
and rearrange to obtain
\begin{equation}\label{eq:extension-growth}
 L\ge\frac{1-\delta}{\delta+2K_G/(\delta\sqrt d)}
 =\frac{d^{1/4}-1}{1+2K_G}
 \ge\frac{d^{1/4}}{2+4K_G}.
\end{equation}
\end{proof}

\begin{proof}[Proof of Theorem~\ref{thm:integral}]
Choose $d_j\to\infty$ and apply Lemma~\ref{lem:related-integral-block} in pairwise disjoint
coordinate blocks of $c_0$. Translate both each $M_{d_j}$ and its
image so that a selected point is zero in both spaces. Multiply the
source and image in block $j$ by the same positive scalar so that the
source diameter is at most $2^{-j}$. Denote the resulting finite set
by $X_j$ and the map by $T_j:X_j\to Y_{d_j}$. Translation and this
simultaneous dilation preserve both the summing norm and the least
extension constant. Because each $X_j$ contains the origin and has diameter at most
$2^{-j}$, the union $X=\bigcup_jX_j$ is a compact countable subset
of $c_0$\amend{. A sequence contained in a finite union of blocks
has a constant subsequence. Otherwise it has a subsequence
$x_{k_r}\in X_{j_r}$ with $j_r\to\infty$, and
$\|x_{k_r}\|_\infty\le 2^{-j_r}$ gives $x_{k_r}\to0\in X$.}
Let $Y=(\bigoplus_jY_{d_j})_{\ell_2}$ and define
$T:X\to Y$ by $T|_{X_j}=T_j$ in the $j$th coordinate.
The map is Lipschitz. Indeed, within a block it is a contraction;
if $x\in X_j$ and $y\in X_k$ with $j\ne k$, then
$\|T(x)-T(y)\|_Y\le(\|x\|_\infty^2+\|y\|_\infty^2)^{1/2}
\le\sqrt2\,\|x-y\|_\infty$.
The target is
separable and reflexive. It also has Rademacher cotype two.
Indeed, with $\delta$ and $Q$ as in
\eqref{eq:finite-target-norm}, consider the linear map
\begin{equation}\label{eq:cotype-embedding}
 J_d:Y_d\longrightarrow\ell_1^S\oplus_1L_1(\Omega),\qquad
 J_du=\bigl(\delta u/2,(\langle r,Qu\rangle)_{r\in\Omega}\bigr).
\end{equation}
The first-moment estimate \eqref{eq:rademacher-first-moment},
combined with $\E|\langle r,Qu\rangle|\le\|Qu\|_2$, shows that
\begin{equation}\label{eq:cotype-distortion}
 \frac{1}{\sqrt3}\|u\|_{Y_d}
 \le\|J_du\|\le2\|u\|_{Y_d}.
\end{equation}
Thus $Y_d$ embeds into an $L_1$ space with distortion at most
$2\sqrt3$, uniformly in $d$. For a finite family $v_i$ in any real $L_1$ space,
Minkowski's integral inequality followed by
\eqref{eq:rademacher-first-moment} yields
\begin{equation}\label{eq:lone-cotype}
 \Big(\sum_i\|v_i\|_1^2\Big)^{1/2}
 \le\int\Big(\sum_i|v_i|^2\Big)^{1/2}
 \le\sqrt3\,\E\Big\|\sum_i\eta_i v_i\Big\|_1,
\end{equation}
for independent fair signs $\eta_i$. By \eqref{eq:cotype-distortion} and \eqref{eq:lone-cotype},
for any finite family $(u_i)$ in $Y_d$,
\begin{equation}\label{eq:uniform-cotype}
 \left(\sum_i\|u_i\|_{Y_d}^2\right)^{1/2}
 \le 6\left(\E\left\|\sum_i\eta_i u_i\right\|_{Y_d}^2\right)^{1/2}.
\end{equation}
For a finite family $(v_i)$ in $Y$, apply
\eqref{eq:uniform-cotype} to the $j$th coordinate of each $v_i$,
square, and sum over $j$. Tonelli's theorem gives
\[
 \sum_i\|v_i\|_Y^2
 \le36\sum_j\E\left\|\sum_i\eta_i(v_i)_j\right\|_{Y_{d_j}}^2
 =36\E\left\|\sum_i\eta_i v_i\right\|_Y^2.
\]
Thus $Y$ has Rademacher cotype two with a constant independent
of the dimensions of its blocks.

In each block take an independent copy $g_j$ of the random
contraction from the proof of Lemma~\ref{lem:related-integral-block}, translating and scaling it with the block and subtracting
its value at zero. Then $g_j(0)=0$ and
$\E|g_j(x)-g_j(y)|\ge\|T_j(x)-T_j(y)\|/4$. For independent fair signs $\eta_j$, independent also of the $g_j$, set
$g(x)=\eta_jg_j(x)/2$ on $X_j$. This is well defined at the common
point $0$ because each $g_j(0)=0$. On each block its Lipschitz
constant is at most $1/2$. For
$x\in X_j$ and $y\in X_k$ with $j\ne k$,
\[
 |g(x)-g(y)|\le(\|x\|_\infty+\|y\|_\infty)/2\le\|x-y\|_\infty.
\]
The identity $(|a+b|+|a-b|)/2=\max\{|a|,|b|\}$ gives on such pairs
\[
 \E|g(x)-g(y)|\ge\tfrac14\big(\E|g_j(x)|+\E|g_k(y)|\big)
 \ge\tfrac1{16}\|T(x)-T(y)\|_Y.
\]
Within a block, \eqref{eq:random-contractions} gives the
corresponding lower bound with constant $1/8$. Hence
$\pi_1^{\mathrm L}(T)\le16$. A Lipschitz extension of $T$ to $c_0$
would, by restriction to the $j$th coordinate block and projection
onto $Y_{d_j}$, extend $T_j$ with the same uniform bound for every
$j$. The lemma rules this out. Since $Y$ is reflexive the same holds
for extensions into $Y^{**}$.

Suppose $T$ were Lipschitz $q$-integral for some $1\le q<\infty$ and
write the factorization as in \eqref{eq:integral-factorization},
$\delta_Y\circ T=B\circ i_{\infty,q}\circ A$.
Restricting a functional on $\operatorname{Lip}_0(Y)$
to $Y^*$ sends evaluation at $y$ to the canonical image of $y$ in
$Y^{**}$ and has norm at most one.

To extend $A$, choose representatives $a_x$ for $A(x)$. Since $X$ is
countable, the inequalities $|a_x-a_y|\le\Lip(A)\|x-y\|_\infty$
hold simultaneously off one null set, on which we may reset all
representatives to zero. For $z\in c_0$ put
\[
 \widetilde A(z)(\omega)=\inf_{x\in X}
     \big(a_x(\omega)+\Lip(A)\|z-x\|_\infty\big).
\]
The term $x=0$ bounds the infimum from above, while the pairwise
scalar Lipschitz inequalities bound it from below. In particular,
$\widetilde A(z)$ is finite almost everywhere.
The infimum is measurable because $X$ is countable. It satisfies
\[
 a_0(\omega)-\Lip(A)\|z\|_\infty
 \le\widetilde A(z)(\omega)
 \le a_0(\omega)+\Lip(A)\|z\|_\infty.
\]
The upper bound is the term $x=0$; the lower bound follows from
$a_x\ge a_0-\Lip(A)\|x\|_\infty$ and the triangle inequality.
Thus $\widetilde A(z)\in L_\infty(\mu)$. Pointwise comparison
of the infima for $z,z'$ gives
$\|\widetilde A(z)-\widetilde A(z')\|_\infty
 \le\Lip(A)\|z-z'\|_\infty$.
For $z=x\in X$, the simultaneous scalar Lipschitz inequalities
show that $\widetilde A(x)=A(x)$ in $L_\infty(\mu)$.
Composing
$\widetilde A$ with $i_{\infty,q}$, $B$, and restriction to $Y^*$
gives the forbidden extension of $T$ into $Y^{**}$. Thus $T$ is
Lipschitz $1$-summing and fails to be Lipschitz $q$-integral for
every finite $q$.
\end{proof}

\section{Six points and the conjugate-exponent test}\label{sec:six}

Let $M$ be the vertex set of the complete bipartite graph $K_{3,3}$,
with sides $P$ and $Q$ and with each edge of length one. Distinct
vertices on the same side have distance two. The identity of $M$ has summing norm $3^{1/p}$, while its
precompositions with contractions on subsets of uniformly convex
spaces have a uniformly smaller norm.

\begin{proof}[Proof of Theorem~\ref{thm:six}]
Give weight one to each of the six unordered pairs on the same side.
For a scalar contraction $g$ let $r_P$ and $r_Q$ be the ranges of
its values on the two sides, measured from minimum to maximum.
The two inequalities $\max_P g-\min_Q g\le1$ and
$\max_Q g-\min_P g\le1$ give $r_P+r_Q\le2$. For three ordered
real values of range $r$, their three pairwise differences in
absolute value are $s,r-s,r$ for some $s\in[0,r]$. The sum of
their $p$th powers is at most $2r^p$. Thus
\[
 \sum_{\{u,v\}\subset P\,\text{or}\,Q}|g(u)-g(v)|^p
 \le 2(r_P^p+r_Q^p)\le 2^{p+1}.
\]
For this particular weighted test, the sum of the $p$th powers
of the six distances under the identity equals $6\cdot2^p$.
By \eqref{eq:lippsum},
\begin{equation}\label{eq:k33-lower}
 \pL(\Id_M)^p\ge\frac{6\cdot2^p}{2^{p+1}}=3.
\end{equation}

For the reverse bound, choose one side with probability $1/2$ and
a uniform singleton on that side. Assign $-1$ to the singleton,
$+1$ to its two companions, and $0$ to the opposite side. Each
realization is a scalar contraction. For a fixed pair $u,v$ on the same side, the probability that
the chosen side is theirs and the singleton separates them is
$(1/2)(2/3)=1/3$. Hence
\begin{equation}\label{eq:k33-random}
 \E|g(u)-g(v)|^p=
 \begin{cases}
 2^p/3,&u,v\text{ distinct on the same side},\\
 1,&u,v\text{ on opposite sides}.
 \end{cases}
\end{equation}
In either case $d_M(u,v)^p\le3\E|g(u)-g(v)|^p$.
For a finite weighted test, summing these pointwise inequalities and
using that the supremum over contractions dominates the expectation
gives $\pL(\Id_M)^p\le3$. Together with \eqref{eq:k33-lower},
this proves $\pL(\Id_M)=3^{1/p}$.

Let $V$ be uniformly convex. There is an $\varepsilon>0$ depending only on $V$ (and,
if necessary, decreased depending on $p$) such that every closed ball
of radius $1+\varepsilon$ in $V$ contains at most two points with
mutual distances at least two. Indeed, choose $\varepsilon<1$ with
$(1+\varepsilon)^{-1}>1-\delta_V(1/2)$, where $\delta_V(t)$ is the
infimum of $1-\|(v+w)/2\|$ over unit-ball vectors at distance at least
$t$. If $u_1,u_2,u_3$
are the three points, subtract the center and divide by
$1+\varepsilon$ to obtain vectors $v_j$ in the unit ball. Since
$\|(v_1-v_j)/2\|\ge (1+\varepsilon)^{-1}$ for $j=2,3$, uniform
convexity applied to $v_1,-v_j$ gives $\|v_1+v_j\|<1/2$. This
would imply $\|v_2-v_3\|<1$, although
$\|v_2-v_3\|\ge 2/(1+\varepsilon)>1$. Decrease $\varepsilon$ further so
that $(1-\varepsilon)^p>1/3$.

Fix a finite $A\subset V$ and a contraction $S:A\to M$. Use the random scalar contractions defined in
\eqref{eq:k33-random}. On the inverse image of the
chosen side use the original values $-1,+1,+1$. Let $D$ be the
inverse image of the singleton and choose, independently for each
label $j$ on the other side, a fair bit $\beta_j\in\{0,1\}$.
At a point $z$ with $S(z)=j$ put
\[
 g(z)=\beta_j\min\{\varepsilon,(\operatorname{dist}(z,D)-1)_+\},
\]
using the value $\beta_j\varepsilon$ if $D$ is empty. Within each
fiber this function is a contraction, and between distinct labels
on the second side its oscillation is at most $\varepsilon$.
Across the two sides its increment is at most one if the first
value is $+1$. If that value is $-1$, the distance to $D$ bounds
the increment. Thus $g$ is a scalar contraction on $A$.

Consider $z,z'$ with distinct labels on the same side. When their
side is chosen, the contribution to $\E|g(z)-g(z')|^p$, including
the probability of this choice, is $2^p/3$. When their side is the
other one, at least one of the three possible singleton fibers
misses the closed ball $B(z,1+\varepsilon)$. Otherwise that ball
would contain three points separated by at least two. With
conditional probability at least $(1/3)(1/4)$ this fiber is
chosen and $\beta_{S(z)}=1$, $\beta_{S(z')}=0$. The event just described has probability at least
$(1/2)(1/3)(1/4)=1/24$, and on that event
$|g(z)-g(z')|=\varepsilon$. Therefore,
\begin{equation}\label{eq:uc-same-side}
 \E|g(z)-g(z')|^p\ge\frac{2^p}{3}+\frac{\varepsilon^p}{24}.
\end{equation}
Across the two sides every increment is at least $1-\varepsilon$.
Consequently, setting
\begin{equation}\label{eq:uc-constant}
 c_{p,V}=\min\left\{\frac13+\frac{\varepsilon^p}{24\cdot2^p},
                      (1-\varepsilon)^p\right\}>\frac13,
\end{equation}
we have $\E|g(z)-g(z')|^p\ge c_{p,V} d_M(Sz,Sz')^p$ for all pairs.
Summation over arbitrary weighted tests gives $\pL(S)\le c_{p,V}^{-1/p}$,
uniformly in the finite set $A$ and the contraction $S$.
Since $c_{p,V}>1/3$, we have
$c_{p,V}^{-1/p}<3^{1/p}=\pL(\Id_M)$.
For $1<p<\infty$, the space $\ell_{p'}$ is uniformly convex
\cite{Clarkson}, and the stated failure of the testing formula follows.
\end{proof}

\begin{corollary}\label{cor:two-triangles}
For $p=2$ the conjugate-exponent test fails even for a contraction
from $M$ to the real Hilbert space $\R^4$.
\end{corollary}

\begin{proof}
Map the two sides to centered equilateral triangles of circumradius
$1/\sqrt2$ in orthogonal planes. The resulting map $T$ has
same-side distances $\sqrt{3/2}$ and cross distances one. The six
same-side pairs give $\pi_2^{\mathrm L}(T)^2\ge9/8$. The first
random scalar test in the proof gives the reverse inequality, since
$(3/2)/(4/3)=9/8$ and the cross distances are one. For the second
test take $\varepsilon$ small enough that $(1-\varepsilon)^2>8/9$.
Every eligible precomposition then satisfies
\[
 \pi_2^{\mathrm L}(T\circ S)^2\le
 \max\left\{\frac{3/2}{4/3+\varepsilon^2/24},
                 \frac1{(1-\varepsilon)^2}\right\}<\frac98.
\]
This gives the strict inequality with $\pi_2^{\mathrm L}(T)=3/(2\sqrt2)$.
\end{proof}

\section{Exact testing on finite subsets of \texorpdfstring{$\ell_1$}{l1}}\label{sec:one}

\begin{theorem}\label{thm:one}
For every Lipschitz map $T:X\to Y$ between metric spaces and every
$1\le p<\infty$,
\[
 \pL(T)=\sup\{\pL(T\circ S): A\subset\ell_1\text{ finite},\quad
                              S:A\to X,\ \Lip(S)\le1\}.
\]
\end{theorem}

\begin{proof}
Write $K$ for the supremum. Precomposition with a contraction gives
$K\le\pL(T)$, so we may assume $K<\infty$. Fix a finite set
$E\subset X$ with at least two points. Adjoin the distinguished point
of $X$ if basepoint preservation is required. Write $o$ for this point,
or choose any $o\in E$ in the unpointed case.

Let $J$ be the edge set of the complete graph on $E$, with edge
weights $w_{\{x,y\}}=d_X(x,y)$. Identify $W=\{0,1\}^J$ with the
family of edge sets and give it the weighted Hamming distance
\[
 d_W(u,v)=\sum_{j\in u\triangle v}w_j.
\]
Here $u\triangle v$ denotes symmetric difference. Define
$\jmath:W\to\ell_1^J$ by
$\jmath(u)=(w_j\mathbf1_{\{j\in u\}})_{j\in J}$.
Then
\begin{equation}\label{eq:hamming-isometry}
 \|\jmath(u)-\jmath(v)\|_1=d_W(u,v),
\end{equation}
so $W$ is a finite metric subspace of $\ell_1^J$. Write $\partial u$ for the set of
vertices of odd degree in $u$. For $x\ne o$ take
$H_x=\{u:\partial u=\{o,x\}\}$, and take
$H_o=\{u:\partial u=\varnothing\}$, which contains the empty edge
set. For distinct $x,y$, toggling their edge is a bijection from
$H_x$ to $H_y$ with displacement exactly $d_X(x,y)$, since it
changes the degree parity at just $x$ and $y$. In particular, all
fibers are nonempty. This is the coset matching of Khot and Naor
\cite[Lemma~3.6]{KN}, here for the kernel of the parity map.

For arbitrary $u\in H_x$ and $v\in H_y$ with $x\ne y$, the
only odd-degree vertices of $u\triangle v$ are $x$ and $y$.
Every connected component of a finite graph has an even number of
odd-degree vertices. Hence $x$ and $y$ belong to the same connected
component of $u\triangle v$, which contains a path joining them.
The sum of the edge weights along that path is at least $d_X(x,y)$
by the triangle inequality. In particular,
\begin{equation}\label{eq:fiber-separation-l1}
 d_W(u,v)\ge d_X(x,y)
 \qquad(u\in H_x,\ v\in H_y,\ x\ne y).
\end{equation}

Let $\mu_x$ be uniform probability on $H_x$ and take
$A_m=\bigcup_{x\in E}H_x^m$ inside $W^m$ with the averaged sum metric.
Concatenating the $m$ blocks with weights $1/m$ embeds $A_m$ in
$\ell_1^{m|J|}$. Sending $H_x^m$ to $x$ defines a contraction $S_m:A_m\to X$
by the lower distance bound, and $S_m(0)=o$. Let $g:A_m\to\R$ be any scalar contraction, and define
$G:E\to\R$ by $G(x)=\int_{H_x^m}g\,d\mu_x^{\otimes m}$.
The edge-toggle matching $H_x\to H_y$ is measure preserving and
moves each point by precisely $d_X(x,y)$. Apply its $m$-fold
product to obtain a coupling whose paired points in $H_x^m$ and
$H_y^m$ have averaged distance $d_X(x,y)$. It follows that
\begin{equation}\label{eq:lone-mean}
 |G(x)-G(y)|\le d_X(x,y)\qquad(x,y\in E).
\end{equation}
Subtracting $G(o)$ yields a pointed scalar test without changing
any increment.

Put $B=\operatorname{diam}W$. On $H_x^m$, changing one coordinate
changes $g$ by at most $B/m$. Let $D_k$ be the difference between the $k$th and
$(k-1)$st conditional expectations in the Doob martingale obtained
by revealing the first $k$ coordinates. Then $|D_k|\le B/m$.
The martingale differences are orthogonal in $L_2$, whence
\begin{equation}\label{eq:lone-variance}
 \int_{H_x^m}|g-G(x)|^2\,d\mu_x^{\otimes m}
 =\sum_{k=1}^m\int|D_k|^2\,d\mu_x^{\otimes m}
 \le\frac{B^2}{m}.
\end{equation}
Since $|g-G(x)|\le B$, for $1\le p\le2$ the $L_p$ norm
is dominated by the $L_2$ norm; for $p\ge2$ we use
$|g-G(x)|^p\le B^{p-2}|g-G(x)|^2$. In both cases,
\begin{equation}\label{eq:lone-oscillation}
 \|g-G(x)\|_{L_p(\mu_x^{\otimes m})}
 \le Bm^{-1/\max\{2,p\}}=:\eta_m.
\end{equation}

Fix pairs $(x_i,y_i)$ in $E$ with weights $a_i\ge0$, and denote
their scalar supremum by $D_E$. Lift each pair by drawing independently
from $\mu_{x_i}^{\otimes m}$ and $\mu_{y_i}^{\otimes m}$, multiplying
the product probabilities by $a_i$. The target distances of $T\circ S_m$
are exactly the original ones. By \eqref{eq:lone-mean}, the mean function $G$ is a scalar
contraction on $E$, so its increments contribute at most $D_E$.
Subtracting the two means and applying Minkowski's inequality
and \eqref{eq:lone-oscillation} gives
\begin{equation}\label{eq:lone-lifted-test}
 \left(\sum_i a_i\,d_Y(Tx_i,Ty_i)^p\right)^{1/p}
 \le K\left(D_E+2\eta_m\Big(\sum_i a_i\Big)^{1/p}\right).
\end{equation}
Let $m\to\infty$. Real scalar contractions on $E$ extend to $X$
with the same constant, so the resulting inequality gives $\pL(T)\le K$.
\end{proof}

\begin{corollary}\label{cor:loneglobal}
If $X$ is Banach and $T:X\to Y$ is Lipschitz, then
\[
 \pL(T)=\sup_{n\ge1}\ \sup_{\substack{F:\ell_1^n\to X\\ \Lip(F)\le1}}
                \pL(T\circ F)
       =\max_{\substack{F:\ell_1\to X\\ \Lip(F)\le1}}\pL(T\circ F).
\]
The maximum is attained also when $\pL(T)=\infty$.
\end{corollary}

\begin{proof}
By Theorem~\ref{thm:testing} it suffices to recover the norm of each
$T\circ U$ for a pointed contraction $U:\ell_\infty^n\to X$.
Apply Theorem~\ref{thm:one} to $T\circ U$. Its proof uses pointed
finite sets $A\subset\ell_1^m$ and pointed contractions
$S:A\to\ell_\infty^n$. Extend each coordinate of $S$ from $A$ to
$\ell_1^m$ with Lipschitz constant at most one. The composition of
this extension with $U$ is a contraction on the whole of $\ell_1^m$
and agrees with $U\circ S$ on $A$. This proves the first equality.

Choose pointed contractions $F_j:\ell_1^{n_j}\to X$ whose
composed norms approach the supremum, including when the supremum is
infinite. Such contractions can be taken to preserve zero: in the construction of Theorem~\ref{thm:one} the distinguished
source point is sent to the distinguished target point, and coordinatewise scalar Lipschitz
extension preserves this property. Identify
$\ell_1=(\bigoplus_j\ell_1^{n_j})_{\ell_1}$ isometrically and put
\begin{equation}\label{eq:lone-direct-sum}
 F((z_j)_j)=\sum_j F_j(z_j).
\end{equation}
The series converges absolutely since $\|F_j(z_j)\|\le\|z_j\|_1$.
The same estimate on differences gives $\Lip(F)\le1$, and restriction
to the $j$th coordinate block recovers $F_j$. Thus $\pL(T\circ F)$
attains the supremum.
\end{proof}

\amend{Corollary~\ref{cor:loneglobal} also yields the Hamming-cube
formula of Theorem~\ref{lone:thm:recovery}. Fix a bounded linear
$U:X\to Y$ and a contraction $F:\ell_1^n\to X$. On a grid
$a+\delta k$, $0\le k_j\le m_j$, encode each integer $k_j$ by
$k_j$ leading ones followed by zeros in a block $B_j$ of $m_j$ bits.
This multiplies distances by $\delta^{-1}$. The map
$D(z)=a+\delta\big(\sum_{b\in B_j}z_b\big)_{j=1}^n$ recovers the
grid points, and $z\mapsto F(D(z))/\delta$ is contractive on the
whole cube. On an encoded test, linearity of $U$ scales the numerator
by $\delta^{-1}$, and scalar Lipschitz extension gives the same factor
in the denominator. Arbitrary tests follow by approximation. Moving
each endpoint by at most $\varepsilon$ changes the scalar supremum by
at most $2\varepsilon(\sum_i w_i)^{1/p}$ and the numerator by at most
$\|U\|$ times this amount. Thus the supremum of the cube tests bounds
$\pL(U\circ F)$ for every $F$ and $n$, and hence bounds
$\pL(U)=\pi_p(U)$ by Corollary~\ref{cor:loneglobal}
and \cite[Theorem~2]{FJ}. Precomposition gives the reverse inequality.}

\appendix
\section{A sharp example on the three-dimensional Hamming cube}\label{lone:sec:three}

Saleh \cite{Saleh} gave a linear program for Lipschitz summing norms
on finite metric spaces. On the Hamming cube the scalar functions
in its constraints may all be taken integer valued. To see this,
let $\mathcal P_N$ be the set of unordered pairs of distinct vertices
and let $\mathcal G_N$ consist of the integer-valued contractions on
$Q_N$ vanishing at the zero vertex. Write $\Delta_e F=F(u)-F(v)$
for $e=\{u,v\}$ in either orientation. For $1\le p<\infty$ the
number $\pL(F)^p$ is the maximum of
\[
 \sum_{e\in\mathcal P_N}w_e\|\Delta_e F\|^p
 \quad\text{over }w_e\ge0\text{ such that }\quad
 \sum_{e\in\mathcal P_N}w_e|\Delta_e g|^p\le1
 \quad(g\in\mathcal G_N).
\]
Indeed, after normalizing the scalar supremum to one, the definition
gives this program with all real scalar contractions vanishing at
zero. These form a compact polytope cut out by the edge inequalities.
At any extreme point, the edges on which the increment has absolute
value one form a connected spanning graph. Otherwise, adding or
subtracting a small constant on a component missing zero preserves
every edge inequality, contradicting extremality. Summing increments
along paths from zero now shows that each extreme point is integer
valued. Since $g\mapsto\sum_e w_e|\Delta_e g|^p$ is convex, its
maximum occurs at an extreme point. Finally, the constraints bound
$w_{\{u,v\}}\le \dH(u,v)^{-p}$, as one sees by testing against the
distance to $u$ minus its value at zero. The feasible set is compact, so the linear program attains its
maximum.

Four pairs distinguish an unrestricted scalar test from the
translation-invariant tests of Proposition~\ref{lone:prop:invariant}.
On $Q_3=\{0,1\}^3$ take
\[
 \mathcal M=\big\{\{000,110\},\{011,101\},\{001,010\},\{100,111\}\big\},
 \qquad
 A(F)=\frac14\sum_{\{u,v\}\in\mathcal M}\|F(u)-F(v)\|.
\]

\begin{proposition}\label{lone:prop:three}
The supremum of $A$ over scalar contractions equals one. Its
supremum over contractions into real Hilbert spaces is
$\sqrt{1+\sqrt2}$, attained in $\R^3$.
\end{proposition}

\begin{proof}
For a scalar contraction, label the greater endpoint in each pair
positive and the lesser endpoint negative, breaking ties arbitrarily.
Each parity class contains two positive and two negative vertices.
Match the two positive even vertices to the two negative odd vertices
by cube edges. This is possible because each even vertex has exactly
one odd nonneighbor, its antipode, and the two antipodes are distinct.
The same argument matches positive odd vertices to negative even
ones. The sum of the four decrements is the sum of the differences
along these four oriented edges, hence at most four. Equality holds
for $g(u)=u_1+u_2+u_3$, whose absolute increments are $2,0,0,2$.

For the Hilbert upper bound write $x_j=(-1)^{u_j}$ and let $M$ be
the involution pairing the endpoints. It flips coordinates $1,2$
on $\{x_1x_2x_3=1\}$ and coordinates $2,3$ on
$\{x_1x_2x_3=-1\}$. Let $E_i(F)$ denote the sum of squared
increments over the four cube edges in direction $i$, and $D(F)$
the sum of the four squared increments over $\mathcal M$. The quadratic estimate
\begin{equation}\label{lone:eq:three-quad}
 D(F)\le\frac{E_1(F)+E_3(F)}{\sqrt2}+E_2(F)
\end{equation}
holds for every Hilbert-valued map $F$.
It suffices to prove this for scalar $f$ and then sum over an
orthonormal basis of the span of the range. With counting measure
on the eight vertices, denote the operators induced by the three
coordinate flips and by $M$ by $U_i$ and $U_M$. For any one of
these involutions $U$, the quadratic form $\langle f,(I-U)f\rangle$
equals the sum of squared increments, each pair counted once.
The asserted inequality is therefore the positivity of
\[
 P=\sqrt2 I-U_2-\frac{U_1+U_3}{\sqrt2}+U_M .
\]
Fix a Walsh character $\chi_B$, put $s_i=(-1)^{\mathbf1_{\{i\in B\}}}$,
and write $B^c=\{1,2,3\}\setminus B$. The parity character
interchanges $\chi_B$ and $\chi_{B^c}$ by multiplication. Thus
$U_M\chi_B=\alpha\chi_B+\beta\chi_{B^c}$, where
$\alpha=s_2(s_1+s_3)/2$ and $\beta=s_2(s_1-s_3)/2$.
These two coefficients are unchanged on replacing $B$ by $B^c$.
On their two-dimensional span the matrix of $P$ is
\[
 \begin{pmatrix}a-t&\beta\\ \beta&a+t\end{pmatrix},
 \qquad a=\sqrt2+\alpha,\qquad t=s_2+(s_1+s_3)/\sqrt2.
\]
If $s_1=s_3$, then $\beta=0$ and $|t|=\sqrt2+s_1s_2=a$.
If $s_1=-s_3$, then $a=\sqrt2$ and $t^2+\beta^2=2$. In both
cases the eigenvalues are $0$ and $2a>0$. These spans exhaust the
Walsh basis. Thus \eqref{lone:eq:three-quad} holds. For a contraction
each $E_i(F)\le4$, so Cauchy--Schwarz gives
$A(F)^2\le D(F)/4\le1+\sqrt2$.

To attain the bound, choose positive $a,b,c$ with
\[
 a^2=(2+\sqrt2)/8,\qquad b^2=\sqrt2/8,\qquad c^2=(2-\sqrt2)/8
\]
and define
\[
 F(x)=
 \begin{cases}
  (a x_1,b x_2,c x_3),&x_1x_2x_3=1,\\
  (c x_1,b x_2,a x_3),&x_1x_2x_3=-1 .
 \end{cases}
\]
On edges in directions $1$ and $3$, the squared length is
$(a+c)^2+(a-c)^2=2(a^2+c^2)=1$. In direction $2$ it is
$2(a-c)^2+4b^2=1$. The edge bounds imply contractivity for Hamming
distance. Each of the four chosen pairs has length
$2\sqrt{a^2+b^2}=\sqrt{1+\sqrt2}$.
\end{proof}

\section*{Acknowledgments}
D. Pellegrino is supported by Grants No. 406457/2023-9 (CNPq/MCTI No. 10/2023), No. 403964/2024-5 (MCTI/CNPq No. 16/2024), and No. 305807/2025-0 from the Conselho Nacional de Desenvolvimento Cient\'ifico e Tecnol\'ogico (CNPq, Brazil). E. Teixeira gratefully acknowledges support from the Grayce B. Kerr Chair at Oklahoma State University and KASUT research grant PRFS-2024-CRG12-6430.2. This work was conducted in part within the DARPA ExpMath project \emph{A Human-Centered Framework for AI-Mathematician Collaboration in Research-Level Mathematics} (Agreement No. HR0011262E029), in which E. Teixeira serves as co-PI. The views and conclusions expressed here are those of the authors and should not be interpreted as representing the official policies of the Department of Defense or the U.S. Government.

\section*{Declaration on generative AI assistance}
\amend{The authors initiated this study of Farmer and Johnson's
Problems~2, 4, and~5 on Lipschitz $p$-summing maps. For the two
Grothendieck questions, they used ChatGPT Astra~6 (OpenAI) to develop
and test Hamming-cube constructions for maps on $\ell_1$ and graph
lifts for the $C(K)$ argument. They also used the tool to check the
scalar Lipschitz estimates on $K_{3,3}$ in the counterexample to the
$\ell_{p'}$ test and to examine how the trace estimate forces extension
constants to grow in the example separating summability from integrality.
Further assistance included locating references, drafting and revising
portions of the exposition, and routine \LaTeX{} editing. The authors
made all final mathematical and editorial decisions and take full
responsibility for the results, proofs, calculations, and references.}

\end{document}